\documentclass[aap,preprint]{imsart}

\RequirePackage{amsthm,amsmath,amsfonts,amssymb}
\RequirePackage[numbers]{natbib}
\RequirePackage{graphicx}

\startlocaldefs
\theoremstyle{plain}
\newtheorem{theorem}{Theorem}[section]
\newtheorem{lemma}[theorem]{Lemma}
 \newtheorem{corollary}[theorem]{Corollary}
 \newtheorem{proposition}[theorem]{Proposition}
 \newtheorem{example}[theorem]{Example}
 \newtheorem{Definition}[theorem]{Definition}
 \theoremstyle{remark}
 \newtheorem{remark}[theorem]{Remark}

\theoremstyle{remark}

\endlocaldefs

\usepackage{bbm}
\usepackage{mathrsfs}

\usepackage{CJK}
\usepackage{amsfonts,amssymb,amsmath,mathrsfs,multirow,booktabs}

 \def\eqref#1{{\rm(\ref{#1})}}

\def\({\left(}\def\){\right)}

\newcommand{\ddr}{\mathrm{d}}

\begin{document}
\begin{frontmatter}
\title{A phase transition in the coming down from infinity of simple exchangeable fragmentation-coagulation processes}
\runtitle{A phase transition in simple EFC processes}

\begin{aug}
\author[A]{\snm{Cl\'ement Foucart}\ead[label=e1]{foucart@math.univ-paris13.fr}}
\address[A]{LAGA, Institut Galil\'ee,
Universit\'e Sorbonne Paris Nord,
\printead{e1}}
\end{aug}

\begin{abstract}
We consider the class of exchangeable fragmentation-coagulation (EFC) processes where coagulations are multiple and not simultaneous, as in a $\Lambda$-coalescent, and fragmentation dislocates at finite rate an individual block into sub-blocks of infinite size. We call these partition-valued processes, simple EFC processes, and study the question whether such a process, when started with infinitely many blocks, can visit partitions with a finite number of blocks or not. When this occurs, one says that the process comes down from infinity.  We introduce two sharp parameters $\theta_{\star}\leq \theta^{\star} \in [0,\infty]$, so that if 
$\theta^{\star}<1$, the process comes down from infinity and if $\theta_{\star}>1$, then it stays infinite. We illustrate our result with regularly varying coagulation and fragmentation measures. In this case, the parameters $\theta^{\star}$, $\theta_\star$ coincide and are explicit. 
\end{abstract}

\begin{keyword}[class=MSC2010]
\kwd[Primary ]{60J25, 60J50}
\kwd{60J80} \kwd{60J90}
\kwd[; secondary ]{60G09}
\end{keyword}


\begin{keyword}
\kwd{Coalescence}
\kwd{Fragmentation}
\kwd{Coming down from infinity}
\kwd{Boundary behavior}
\kwd{Partition-valued processes}
\kwd{Exchangeability}
\kwd{Coupling}
\end{keyword}

\end{frontmatter}

\section{Introduction and main results}

Fragmentation and coagulation are natural phenomena that can be observed in many different contexts. We refer to Bertoin \cite{coursbertoin} and Pitman \cite{coursPitman} for an introduction to exchangeable fragmentations and coalescents. These processes form random systems of disjoint subsets, so-called blocks, covering $\mathbb{N}:=\{1,2...\}$, evolving either by fragmentations of blocks or by coagulations of two or more blocks.
By exchangeability, it is meant that the rate of coalescence only depends on the number of subsets that are merging and not on their constituent elements. Similarly, blocks fragmentate into sub-blocks independently of each other, with a same rate.

One striking feature of pure exchangeable coalescents lies in the so-called ``coming down from infinity". This phenomenon states that, although started from a partition with infinitely many blocks, the coalescent process reaches a partition with only a finite number of blocks. This phenomenon has received a great deal of attention in the last two decades. 
The most important results in this respect are certainly Schweinsberg's  necessary and sufficient condition,  \cite{Schweinsberg00}, for the coming down from infinity of coalescents with no simultaneous coagulations, and the study of their speed of coming down by Berestycki et al. \cite{Beres10} and Limic and Talarczik \cite{limic2015}.

Most studies have been carried out for processes of pure fragmentation or pure coagulation. However many natural stochastic particle models, ranging from physics to mathematical genetics, evolve in time by both fragmentation and coalescence.  We refer for instance to Aldous's review \cite[Sections 1.4 and 1.5]{MR1673235} for a list of models.  This led Berestycki to define in \cite{Berestycki04} a class of partition-valued processes called exchangeable fragmentation-coagulation (EFC) processes. 
Some examples of EFC processes have been recently studied by Bertoin and Kortchemski \cite{BertoinKortchemski}, Kyprianou et al. \cite{kyprianou2017} and Foutel-Rodier et al. \cite{foutelrodier2019}. It is also noteworthy that EFC processes arise in the background of several mathematical population models  involving interactions, see for instance \cite{Lambert:2005bq}, \cite{Foucartimpact} and \cite{FoucartEJP}, as well as Gonz\'alez-Casanova and Span\`o \cite{GonzalesSpano}, Gonz\'alez-Casanova et al. \cite{Gonzalesetal} and Foucart and Zhou \cite{FoucartZhouWF}.
\\

The purpose of this article is to consider the coming down from infinity phenomenon for EFC processes. We stress that in the literature the terminology ``coming down from infinity" has been used in different contexts and often includes  the assumption that the boundary $\infty$ is inaccessible for the process under study. We do not assume this here and when an EFC process comes down from infinity, it may also return to a partition with infinitely many blocks at some other positive time.

In his seminal paper, Berestycki has shown that EFC processes are characterized in law by two exchangeable $\sigma$-finite measures, $\mu_{\mathrm{Frag}}$ and $\mu_{\mathrm{Coag}}$ on $\mathcal{P}_\infty$, the space of partitions of $\mathbb{N}$, governing respectively the fragmentation and the coagulation in the system. Among other results, he established in \cite[Theorem 12]{Berestycki04}, that when the fragmentation occurs at infinite rate, namely $\mu_{\mathrm{Frag}}(\mathcal{P}_\infty)=\infty$, the EFC process may have finitely many blocks only at times of exceptional coalescence in which, instantaneously, infinitely many blocks are merged into a finite number. It leads naturally to discard these cases for a further study of the block-counting process. 
In this direction, Kyprianou et al. \cite{kyprianou2017} have studied a particular extreme example, the so-called ``fast" EFC process,  where pairwise coagulations occur at rate $c_{\mathrm{k}}>0$, as in the Kingman coalescent, and fragmentation splits any individual block into its constituent elements at finite rate $\lambda$, creating thus infinitely many singletons blocks. 
They establish a nice phase transition phenomenon, see \cite[Theorem 1.1]{kyprianou2017} stating that the ``fast"-EFC process comes down from infinity if and only if $\theta:=\frac{2\lambda}{c_{\mathrm{k}}}<1$.

We will investigate such properties for a class of EFC processes with less extreme fragmentation and coagulation mechanisms. We call them \textit{simple} EFC processes  and describe them now briefly.  
%
We assume that the fragmentation measure is finite, i.e. $\mu_{\mathrm{Frag}}(\mathcal{P}_\infty)<\infty$, and is supported by the partitions with no singleton blocks.
As we shall notice in the sequel, see Section \ref{EFCdef}, since the measure $\mu_{\mathrm {Frag}} $ is exchangeable, this latter condition is equivalent to the fact that the measure $\mu_{\mathrm {Frag}} $ has for support the set of partitions whose blocks are infinite, that is to say
\begin{equation}\label{infiniteblock}
\mu_{\mathrm{Frag}}(\{\pi; \#\pi_i<\infty \text{ for some }i\leq \#\pi\})=0,
\end{equation} where we have denoted by $\#\pi_i$ the number of elements in the block $\pi_i$ and $\#\pi$ the number of blocks in the partition $\pi$. We assume furthermore that there are no simultaneous multiple coagulations of blocks, nor coagulations of all blocks at once.  Under this latter assumption, coalescences occur as in a $\Lambda$-coalescent. The measure $\Lambda$, governing coalescences, stands for a finite Borel measure on $[0,1)$ of the form \[\Lambda(\ddr x):=x^{2}\nu_{\mathrm{Coag}}(\ddr x)+c_{\mathrm{k}}\delta_{0}\]
where $c_{\mathrm{k}}\geq 0$ is the Kingman parameter, driving pairwise coagulations, and $\nu_{\mathrm{Coag}}$ is a Borel measure \footnote{We shall also call coagulation measure, the measure $\nu_{\mathrm{Coag}}$, but this should not cause any confusion.} on $(0,1)$, driving multiple coagulations and satisfying $\int_{0}^{1}x^{2}\nu_{\mathrm{Coag}}(\ddr x)<\infty$.

We shall establish in the forthcoming Section \ref{simpleEFC} that if $(\Pi(t), t\geq 0)$ is a simple EFC process, then its number of blocks $(\#\Pi(t),t\geq 0)$ is a Markov process when evolving in $\mathbb{N}$, with the following transitions:
\begin{itemize}
\item[(i)] from $n$ to $n+k$, with $k\in \mathbb{N}\cup \{\infty\}$, at rate $n\mu(k)$ where 
\begin{equation*} \mu(k):=\mu_{\mathrm{Frag}}(\{\pi,\ \#\pi=k+1\}),
\end{equation*}
\item[(ii)] from $n$ to $n-k+1$, with $2\leq k\leq n$, at rate $\binom{n}{k}\lambda_{n,k}$ where
\begin{equation*}  \lambda_{n,k}:=\int_{0}^{1}x^{k-2}(1-x)^{n-k}\Lambda(\ddr x)=c_{\mathrm{k}}\mathbbm{1}_{\{k=2\}}+\int_{0}^{1}x^{k}(1-x)^{n-k}\nu_{\mathrm{Coag}}(\ddr x).
\end{equation*}
\end{itemize}
The measure $\mu$ in (i), called splitting measure, is by definition the image of $\mu_{\mathrm{Frag}}$ by the map $\pi \mapsto \#\pi-1$. Note that $\mu(\infty)$ can be positive. A simple use of the exchangeability property, see the forthcoming background section, ensures that $\mu$ can be \textit{any} finite measure on $\mathbb{N}\cup \{\infty\}$. 

Some simple EFC processes have already been studied in the literature. 
When there are no multiple coagulations (namely $\nu_{\mathrm{Coag}}\equiv 0$) but only binary coalescences at rate $c_{\mathrm{k}}>0$, and under the additional assumption that $\mu(\infty)=0$, Berestycki \cite[Section 5]{Berestycki04} and Lambert \cite[Section 2.3]{Lambert:2005bq} have observed that the process $(\#\Pi(t),t\geq 0)$ has the same transitions as a discrete \textit{logistic} branching process (defined in Section 2 of \cite{Lambert:2005bq}), when $\Pi$ starts from a partition with blocks of infinite size. A sufficient condition over $\mu$  (entailing $\mu(\infty)=0$) for coming down from infinity of the EFC process, see \cite[Proposition 15]{Berestycki04}, was derived from this observation. We also wish to mention that continuous-time Markov chains with jump rates (i) and (ii) have been studied in \cite{Gonzalesetal}, under some assumptions on $\Lambda$ and $\mu$.
   
Our main aim is to study the coming down from infinity for the whole class of simple EFC processes. In particular, we shall not make any assumption on the measure $\mu$.
We will find sharp parameters measuring, in some sense, how fragmentation interplays with the coagulations and obtain a general phase transition phenomenon for coming down from infinity. 

Plainly if the pure $\Lambda$-coalescent process stays infinite, then any EFC process with coalescences driven by $\Lambda$ stays infinite. We work therefore, without loss of generality, under the assumption that the pure coalescent comes down from infinity.  Recall  Schweinsberg's condition. For any $n\geq 2$, set 
\begin{equation} \label{phi} \Phi(n):=\sum_{k=2}^{n}\binom{n}{k}\lambda_{n,k}(k-1). \end{equation} 
This is the rate at which the number of blocks is decreasing when the pure $\Lambda$-coalescent process starts from $n$ blocks. The pure $\Lambda$-coalescent comes down from infinity if and only if
\begin{equation}\label{CDI} \sum_{n\geq 2}\frac{1}{\Phi(n)}<\infty \qquad \text{ (Schweinsberg's condition)}. 
\end{equation} 
Fundamental properties of $\Lambda$-coalescents and of the function $\Phi$ are recalled in Section \ref{cdicoal}.\\

Recall the definition of the splitting measure $\mu$ and  for any $k\geq 1$, let $\bar{\mu}(k)$ be its tail $\bar{\mu}(k):=\mu(\{k,k+1,\cdots,\infty\})$. 
\begin{theorem}\label{mainthm} Let $(\Pi(t),t\geq 0)$ be a simple EFC process started from an exchangeable partition such that $\#\Pi(0)=\infty$. Assume \eqref{CDI} and set 
\begin{equation}\label{theta}
\theta_{\star}:=\underset{n\rightarrow \infty}\liminf\sum_{k=1}^{\infty}\frac{n\bar{\mu}(k)}{\Phi(k+n)}\in [0,\infty] \text{ and } \theta^{\star}:=\underset{n\rightarrow \infty}\limsup\sum_{k=1}^{\infty}\frac{n\bar{\mu}(k)}{\Phi(k+n)}\in [0,\infty]. \end{equation}
If 
$\theta^{\star}<1$ then the process comes down from infinity. If $\theta_{\star}>1$ then the process stays infinite.
\end{theorem}

Note that  if $\theta^{\star}=0$, the process comes down from infinity and if $\theta_{\star}=\infty$, it stays infinite. When both parameters agree, namely $\theta^{\star}=\theta_\star$, we shall denote  their common value simply by $\theta$. A phase transition will occur at $\theta$, between the regime where the process stays infinite and the regime where it visits partitions with a finite number of blocks. The parameters $\theta_{\star}$ and $\theta^{\star}$ measure how fragmentations and coagulations combine when there are a large number of blocks. Indeed, $n\bar{\mu}(k)$ is the rate at which the process, started from a partition with $n$ blocks, jumps to a partition with more than $n+k$ blocks and $\Phi(n+k)$ is the decrease rate when there are $n+k$ blocks. A more precise heuristics of Theorem \ref{mainthm} is given in Section \ref{proof}.

The question whether or not the EFC process can reach partitions with infinitely many blocks when it starts from a finite partition is not addressed in this work. Clearly $\mu(\infty)>0$ is a sufficient condition for $\infty$ to be accessible for the process $(\#\Pi(t),t\geq 0)$. The case $\mu(\infty)=0$ is more involved and is studied in Foucart and Zhou \cite{explosion}.\\

The proof of Theorem \ref{mainthm} is based on two different couplings of the partition-valued process $(\Pi(t),t\geq 0)$ and is differed in Section \ref{proof}.\\

As a first application of Theorem \ref{mainthm}, we study the case where the fragmentation can split blocks into infinitely many sub-blocks. 
\begin{corollary}\label{FEFC2} Assume that the measure $\Lambda$ satisfies \eqref{CDI}. Recall $c_{\mathrm{k}}=\Lambda(\{0\})\geq 0$ and set $\lambda:=\mu(\infty) \geq 0$.
\begin{itemize}
\item[(1)]  If $c_{\mathrm{k}}>0$ then $\theta=2\lambda/c_{\mathrm{k}}$. In particular, if $\lambda=0$ then  $\theta=0$.
\item[(2)]  If $\lambda>0$ and $c_{\mathrm{k}}=0$, then $\theta=\infty$. 
\end{itemize}
\end{corollary}
The corollary above ensures that when there are binary coagulations, namely $c_{\mathrm{k}}>0$, a fragmentation measure with no mass on the partitions with infinitely many blocks, i.e. $\mu(\infty)=0$, will never prevent the EFC process to come down from infinity. Moreover when $\mu(\infty)>0$, only coalescences with a Kingman component can make the process come down from infinity. 
\begin{remark}
The phase transition in (1), occurring at $\theta=2\lambda/c_{\mathrm{k}}$, is similar as the one observed in the ``fast"-EFC process in \cite[Theorem 1.1]{kyprianou2017}. We shall see later in Proposition \ref{propcritical}, that when only binary coagulations are allowed, namely $\Lambda=c_{\mathrm{k}}\delta_0$, the process stays infinite in the case $\theta=1$.  
\end{remark}
\begin{corollary}\label{suffcond1}
Assume $\mu(\infty)=0$ and that the measure $\Lambda$ satisfies \eqref{CDI}. 
\begin{equation}\label{(2)}
\text{If }\ \sum_{k=2}^{\infty}\frac{k}{\Phi(k)}\bar{\mu}(k)<\infty, \text{ then } \theta=0.
\end{equation}
\end{corollary}
\begin{remark}
The series convergence  in Corollary \ref{suffcond1} is a tractable sufficient condition, however it is far from being necessary. For instance, when $\Lambda=c_{\mathrm{k}}\delta_0$, $\Phi(k)=c_{\mathrm{k}} \binom{k}{2}$ for all $k\geq 2$ and one can check that the condition of convergence of the series in \eqref{(2)} coincides with a log-moment condition on $\mu$. We know however by Corollary \ref{FEFC2} that when $c_{\mathrm{k}}>0$, such a moment assumption is not necessary in order to have $\theta=0$.  Note also that $(k/\Phi(k),k\geq 2)$ is always bounded, so that if $\mu$ admits a first moment, then the condition in \eqref{(2)} is fulfilled and the process comes down from infinity as soon as \eqref{CDI} holds. We mention that Gonz\'ales et al. 
\cite[Theorem 2-(i)]{Gonzalesetal} have shown that if $\mu$ admits a first moment and $\frac{c_{\mathrm{k}}}{2}>\sum_{k=1}^{\infty}k\mu(k)$ then, the minimal process with jumps (i) and (ii) has $\infty$ as entrance boundary.
\end{remark}
The parameters $\theta^{\star}$ and $\theta_{\star}$, in their very definition \eqref{theta},  are rather intricate. We will provide in Section \ref{proofsandexamples} sufficient conditions entailing either $\theta^{\star}=0$, $\theta_{\star}=\infty$ or $\theta_{\star}$ and  $\theta^{\star}\in (0,\infty)$, see the forthcoming Lemma \ref{boundsigma}. This enables us in particular to find classes of EFC processes with $\theta_{\star}=\theta^{\star}=\theta\in (0,\infty)$. 


\begin{proposition}\label{regularcdi} Let $d>0$ and $\lambda>0$. Assume $\Phi(n)\underset{n\rightarrow \infty}{\sim} d n^{1+\beta}$ and  $\bar{\mu}(n)\underset{n\rightarrow \infty}{\sim}\frac{\lambda}{n^{\alpha}}$, with $\alpha>0$, $\beta\in (0,1)$.
We have the following three cases: 
\begin{enumerate}
\item[(1)] $\beta<1-\alpha$  then $\theta=\infty$ and the process stays infinite, 
\item[(2)]  $\beta>1-\alpha$  then $\theta=0$ and the process comes down from infinity, 
\item[(3)] $\beta=1-\alpha$  then  $\Phi(n)\underset{n\rightarrow \infty}{\sim} d n^{2-\alpha}$ and one has 
\[\theta=\frac{\lambda}{d(1-\alpha)}\in (0,\infty).\]
\end{enumerate}
\end{proposition}
\begin{remark} 
Important examples of coagulation measures satisfying $\Phi(n)\underset{n\rightarrow \infty}{\sim} d n^{1+\beta}$ for some $d>0$ and $\beta\in (0,1)$ are measures $\Lambda$ of the Beta form \begin{equation}\label{beta}\nu_{\mathrm{Coag}}(\ddr x)=x^{-2}\Lambda(\ddr x)=\frac{c}{\text{Beta}(1-\beta,a)}x^{-\beta-2}(1-x)^{a-1}\ddr x
\end{equation}
with $a>0$ and $c>0$. In this case, the factor constant $d$ is $c\frac{\Gamma(a-\beta+1)}{\Gamma(a)\beta(\beta+1)}$ and when $\beta=1-\alpha$, the phase transition occurs at
\[\theta=\frac{(2-\alpha)\lambda}{c} \frac{\Gamma(a)}{\Gamma(a+\alpha)} \in (0,\infty).\]
Heuristically, by letting $\alpha$ towards $0$, one recovers $\theta=\frac{2\lambda}{c}$ the parameter of the phase transition in the case $\mu=\lambda \delta_{\infty}$ and $c_{\mathrm{k}}=c$.  
\end{remark}
We also consider the case of EFC processes with  ``slower" coalescences. 


\begin{proposition}\label{slowcdi} Let $d>0$ and $\lambda>0$. Assume $\Phi(n)\underset{n\rightarrow \infty}{\sim} dn(\log n)^{\beta}$ with $\beta>1$, and $\bar{\mu}(n)\underset{n\rightarrow \infty}{\sim}\lambda \frac{(\log n)^{\alpha} }{n}$ with $\alpha \in \mathbb{R}$. We have the following three cases
\begin{enumerate}
\item[(1)] $\beta<1+\alpha$  then $\theta=\infty$ and the process stays infinite, 
\item[(2)]  $\beta>1+\alpha$  then $\theta=0$ and the process comes down from infinity, 
\item[(3)] $\beta=1+\alpha$  then  $\Phi(n)\underset{n\rightarrow \infty}{\sim} dn\log(n)^{1+\alpha}$ and one has 
\[\theta=\frac{\lambda}{d(1+\alpha)}\in (0,\infty).\]
\end{enumerate}
\end{proposition}
We will explain in Section \ref{cdicoal} how to construct coagulation measures in order to have the equivalences $\Phi(n)\underset{n\rightarrow \infty}{\sim} dn^{\beta+1}$ with $\beta>0$ or $\Phi(n)\underset{n\rightarrow \infty}{\sim} dn(\log n)^{\beta}$ for $\beta\in (1,\infty)$. 
\\


The paper is organized as follows. In Section \ref{background}, we provide some background on exchangeable random partitions, recall the definition of an EFC process, and in particular explain its Poissonian construction. We focus then on simple exchangeable coalescents and simple EFC processes. We show in Section \ref{simpleEFC} that the number of blocks $(\#\Pi(t),t\geq 0)$ has the same dynamics as explained in the introduction. Section \ref{proof} is devoted to the proof of Theorem \ref{mainthm}. The proof is based on several couplings on the space of partitions. We show Corollary \ref{FEFC2}, Corollary \ref{suffcond1}, Proposition \ref{regularcdi}, Proposition \ref{slowcdi} in Section \ref{proofsandexamples}. 

\section{Background on exchangeable fragmentation-coalescence processes}\label{background} 
\subsection{Exchangeable random partitions and EFC processes}\label{EFCdef}
We refer to Bertoin's book \cite[Section 2.3 in Chapter 2]{coursbertoin} for background on exchangeable random partitions. For any $n,m\in \mathbb{N}$ such that $n\leq m$, the integer interval between $n$ and $m$ is denoted by $[|n,m|]$. For any $n\in \mathbb{N}\cup\{\infty\}$, we set $[n]=[|1,n|]$ and call partition of $[n]$, a collection $\pi=\{\pi_1,\pi_2,\cdots\}$ of subsets of $\mathbb{N}$ satisfying 
 $\pi_{i}\cap \pi_{j}=\emptyset \text{ when } i\neq j \text{ and } \cup_{i=1}^{\infty}\pi_i=[n]$.  The blocks of the partition $\pi$ are listed in the order of their least element. Namely, if $\pi_j$ is the $j$-th block of $\pi$, then for any $i\leq j$, $\min \pi_{i}\leq \min \pi_j$. Recall that $\#\pi$ denotes the number of non-empty blocks of $\pi$. By convention, 
if $\#\pi=m<\infty$ then we set $\pi=(\pi_1,\cdots,\pi_m,\emptyset,\cdots)$ where $(\emptyset,\cdots)$ denotes a countably infinite collection of empty sets.  The space of partitions of $[n]$ is denoted by $\mathcal{P}_{n}$. In particular, $\mathcal{P}_\infty$ is the set of partitions of $[\infty]=\mathbb{N}$. Any partition $\pi\in \mathcal{P}_n$ can also be represented as an equivalence relation $\underset{\pi}{\sim}$ over $[n]$ by stating 
\[i\underset{\pi}{\sim} j \text{ if and only if } i,j \text{ belong to the same block of }\pi.\]
For any $m\geq n$ and $\pi\in \mathcal{P}_m$, we denote by $\pi_{|[n]}$ the restricted partition $(\pi_i\cap [n],i\geq 1)$. 
Note that for any partition $\pi$,
$(\#\pi_{|[n]},n\geq 1)$ increases towards $\#\pi$ as $n$ goes to $\infty$. We endow $\mathcal{P}_\infty$, with the compact metric  \begin{equation}\label{distance} d(\pi,\pi')=\left(\max\{n\geq 1, \pi_{|[n]}=\pi'_{|[n]}\}\right)^{-1}.
\end{equation}

For any $n\in \mathbb{N}\cup\{\infty\}$, set $0_{[n]}:=\{\{1\},\cdots, \{n\},\emptyset,\cdots\}$ and $1_{[n]}:=\{[n],\emptyset,\cdots\}$ where we have denoted by $\emptyset,\cdots$  a countable collection of empty sets. We introduce now the operations of coagulation and fragmentation.

\begin{Definition}\label{coagfrag} Let $n\in \mathbb{N}\cup \{\infty\}$ and $m\in 
\mathbb{N}\cup \{\infty\}$, $\pi\in \mathcal{P}_n$ and $\pi'\in \mathcal{P}_m$ and $k\in \mathbb{N}$.
\begin{itemize}
\item If $\#\pi \leq m$, a coagulation of $\pi$ by $\pi'$, denoted by $\mathrm{Coag}(\pi,\pi')$, is  a partition of $[n]$ defined  by
\[\mathrm{Coag}(\pi,\pi'):=\{\underset{j\in \pi'_i}{\cup}\pi_j\ ;\ i\geq 1\}.\]
\item If $\#\pi\geq k$, a fragmentation of the $k$-th block of $\pi$ by $\pi'$, denoted by $\mathrm{Frag}(\pi,\pi',k)$, is the collection of sets 
\[\mathrm{Frag}(\pi,\pi',k):= \left(\{\pi_i\ ; \ i\in [|1,\#\pi|]\setminus \{k\}\}\cup \{\pi_k\cap\pi'_j \ , j\geq 1\}\right)^{\downarrow}\]
where the notation $\left(\dots\right)^{\downarrow}$ means that we are reindexing by their least element the collection of sets formed by the sub-blocks of $\pi_k$ according to $\pi'$ and all $\pi_i$ for $i\neq k$.
\end{itemize}
\end{Definition}
\noindent For instance, let $\pi=\{\{1,3\},\{2,5\}, \{4\}\}$, $\pi'=\{\{1\},\{2,3\}\}$ and $k=1$. Then,
\begin{align*}
\mathrm{Coag}(\pi,\pi')&= \{\{1,3\},\{2,4,5\}\}, \text{ and }\\
\mathrm{Frag}(\pi,\pi',1)&=\{ \{1,3\}\cap \{1\},\{1,3\}\cap\{2,3\},\{2,5\},\{4\}\}^{\downarrow}\\
&=\{\{1\},\{2,5\},\{3\},\{4\}\}.
\end{align*}
For any $\pi\in \mathcal{P}_n, \pi'\in \mathcal{P}_m$ with $m\geq \#\pi$, the partition $\mathrm{Coag}(\pi,\pi')$ is coarser than $\pi$. For any $k\leq \#\pi$, when $m\geq \max \pi_k$, $\pi_k \cap [m]=\pi_k$ and $\cup_{i\geq 1}\mathrm{Frag}(\pi,\pi',k)_i= \cup_{i=1 \atop i\neq k}^{\#\pi} \pi_i \cup \left( \pi_k \cap [m] \right)=[n]$, so that 
$\mathrm{Frag}(\pi,\pi',k)$ is also a partition of $[n]$, which is finer than $\pi$. Lastly, for any partition $\pi$, for which it makes sense, one has $\mathrm{Coag}(\pi,0_{[n]})=\mathrm{Coag}(0_{[n]},\pi)=\pi$  and $\mathrm{Frag}(\pi,1_{[n]},j)=\pi$ for any $j\in [|1, \#\pi|]$. \\

\begin{figure}[h!]
\centering \noindent
\includegraphics[scale=0.9]
{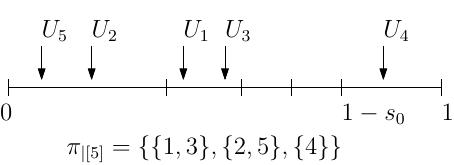}
\caption{paint-box}
\label{paintbox}
\end{figure}

 
Let $\sigma$ be  a permutation of $\mathbb{N}$ with finite support. Namely there is $n\in \mathbb{N}$ such that for any $m\geq n$, $\sigma(m)=m$. The permutation $\sigma$ acts on $\mathcal{P}_\infty$ as follows: we define the partition $\sigma \pi$ by the equivalence relation $i\underset{\sigma \pi}{\sim} j$ if and only if $\sigma(i)\underset{\pi}{\sim} \sigma(j)$. We now recall some elements about exchangeable random partitions. From now on, we shall work on the space $\mathcal{P}_\infty$ equipped with the Borelian $\sigma$-field generated by $d$.
\begin{Definition} A random partition $\pi$ of $\mathbb{N}$ is said to be exchangeable if for any permutation $\sigma$ with finite support the random partitions $\sigma \pi$ and $\pi$ have the same law.
\end{Definition}
A generic example of 
exchangeable random partition is the so-called \textit{paint-box}. 
Define the space of mass-partitions $$\mathcal{P}_{\mathrm{m}}:=\left\{(s_{1},s_{2},...);s_{1}\geq s_{2}\geq ...\geq 0, \sum_{i=1}^{\infty}s_{i}\leq 1\right \}.$$
Let $\mathrm{s}\in \mathcal{P}_{\mathrm{m}}$ and set $s_0=1-\sum_{i=1}^{\infty}s_i$. Partition the interval $[0,1-s_0]$ into subintervals of length $(\mathrm{s}_i,i\geq 1)$. Let  $(U_i,i\geq 1)$ be an i.i.d sequence of uniform random variables over $[0,1]$. The $\mathrm{s}$-paintbox is the random partition $\pi$ defined by letting $i$ and $j$ in the same block if and only if $U_i$ and $U_j$ fall into a same subinterval of $[0,1-s_0]$. When $U_i$ falls into the dust, namely $[1-s_0,1]$, the integer $i$ forms a singleton block of the partition $\pi$ (see Figure \ref{paintbox}).

We denote by $\rho_\mathrm{s}$ the law of the random partition $\pi$. Clearly this random partition is exchangeable, and its law $\rho_{\mathrm{s}}$ does not depend on the locations of the subintervals of $[0,1]$, (including the subinterval $[1-s_0,1]$). Reciprocally, Kingman has shown in \cite{Kingman1} that any exchangeable random partition has the same law as a mixture of paint-boxes. Namely if $\pi$ is a random  exchangeable partition then \[\mathbb{P}(\pi \in \cdot)=\int_{\mathcal{P}_{\mathrm{m}}}\rho_{\mathrm{s}}(\cdot)\nu(\ddr \mathrm{s})\]
where $\nu$ is a probability measure over $\mathcal{P}_{\mathrm{m}}$. The probability measure $\nu$ corresponds to the law of the ranked asymptotic frequencies of $\pi$: $|\pi|^{\downarrow}=(|\pi_i|,i\geq 1)^{\downarrow}$ with
\[|\pi_i|=\underset{n\rightarrow \infty}{\lim}\frac{\#(\pi_{i}\cap [n])}{n} \text{ a.s.}\]
We refer to \cite[Proposition 2.8 and Theorem 2.1 page 100]{coursbertoin} for fundamental properties of random exchangeable partitions. We shall remind us the following generic properties of an $\mathrm{s}$-paint box. If $s_0=0$, then $\pi$ has no singleton block and each block has infinitely many elements. If moreover, $s_i>0$ for all $i\geq 1$, then $\pi$ has infinitely many blocks. Lastly, if there is $k$ such that for all $i\geq k+1$, $s_i=0$, then the partition has at most $k$ blocks. If $s_0>0$, infinitely many random variables $U_i$ will fall in $[1-s_0,1]$ almost surely and there are infinitely many singletons (the so-called dust).

Last, recall that if $\pi$ and $\pi'$ are two independent exchangeable random partitions then $\mathrm{Coag}(\pi,\pi')$ is also exchangeable. Similarly, if one chooses uniformly at random  a block $k$ among those of $\pi$ (in a loose sense, since the partition $\pi$ might have infinitely many blocks) and splits it with $\pi'$, the random partition $\mathrm{Frag}(\pi,\pi',k)$ is exchangeable. This preservation of the exchangeability property naturally leads to consider the following class of processes.

\begin{Definition}[Definition 1 in \cite{Berestycki04}]\label{DefEFC} An exchangeable Fragmentation-Coagulation process is a process $(\Pi(t),t\geq 0)$ valued in $\mathcal{P}_\infty$ satisfying the properties:
\begin{itemize}
\item For any $t\geq 0$, $\Pi(t)$ is exchangeable.
\item For any $n\in \mathbb{N}$, the process $(\Pi_{|[n]}(t),t\geq 0)$ is a c\`adl\`ag Markov chain valued in $\mathcal{P}_n$ evolving by fragmentation of one of its block or by coagulation.
\end{itemize}
More precisely, when the process $\Pi_{|[n]}$ is in state $\pi\in \mathcal{P}_n$, it can only jump to a partition $\pi'$ such that either $\pi'=\mathrm{Coag}(\pi,\pi'')$ for some partition $\pi''$ or  $\pi'=\mathrm{Frag}(\pi,\pi'',k)$ for some partition $\pi''$ and some $k\geq 1$.
\end{Definition} 

It is shown in \cite[Proposition 4]{Berestycki04} that any EFC process is characterized in law by four parameters $(c_{\mathrm{k}}, c_{\mathrm{e}}, \nu_{\mathrm{Coag}}, \nu_{\mathrm{Disl}})$ where $c_{\mathrm{k}}\geq 0, c_{\mathrm{e}}\geq 0$ and $\nu_{\mathrm{Coag}}, \nu_{\mathrm{Disl}}$ are positive $\sigma$-finite measures on $\mathcal{P}_{\mathrm{m}}$, respectively called the \textit{coagulation} and \textit{dislocation} measures. Those measures satisfy the following conditions: $\nu_{\mathrm{Coag}}(\{0,0, \cdots\})=0, \nu_{\mathrm{Disl}}(\{1,0, \cdots\})=0$ and 
\[\int_{\mathcal{P}_\mathrm{m}}\left(\sum_{i\geq 1}s_{i}^{2} \right)\nu_{\mathrm{Coag}}(\ddr s)<\infty \text{ and } \int_{\mathcal{P}_\mathrm{m}}\left(1-\sum_{i\geq 1}s_{i}^{2}\right)\nu_{\mathrm{Disl}}(\ddr s)<\infty.\] The coefficients $c_{\mathrm{k}}$ and $c_{\mathrm{e}}$ are called the Kingman coefficient and erosion coefficient. \\

We now briefly explain the construction of EFC processes. We refer the reader to \cite[Section 3.2]{Berestycki04} for more details.  
\\

\textbf{Poisson construction}. For every pair $(i,j)\in \mathbb{N}^2$ with $i<j$, we write $K(i,j)$ for the partition of $\mathbb{N}$ whose blocks consist of the pair $\{i,j\}$ and the singletons $\{k\}$ for $k\neq i,j$. For any $i\in \mathbb{N}$, let also $e(i)$ be the partition $\{\mathbb{N}\setminus \{i\}, \{i\}\}$. Define the $\sigma$-finite exchangeable measures \begin{equation} 
\label{represmucoag}
\mu_{\mathrm{Coag}}(\ddr \pi):=c_{\mathrm{k}}\sum_{1\leq i<j}\delta_{K(i,j)}+\int_{\mathcal{P}_{\mathrm{m}}}\rho_{s}(\cdot)\nu_{\mathrm{Coag}}(\ddr s)
\end{equation}
and
\begin{equation}\label{represmufrag} \mu_{\mathrm{Frag}}(\ddr \pi):=c_{\mathrm{e}}\sum_{i\geq 1}\delta_{e(i)}+\int_{\mathcal{P}_{\mathrm{m}}}\rho_{s}(\cdot)\nu_{\mathrm{Disl}}(\ddr s).
\end{equation}
Denote by $\#$ the counting measure over $\mathbb{N}$. Consider two independent Poisson point processes $\text{PPP}_{C}=\sum_{t>0}\delta_{(t,\pi^{c})}$ and $\text{PPP}_{F}=\sum_{t>0}\delta_{(t,\pi^{f},k)}$ respectively on $\mathbb{R}_+\times \mathcal{P}_\infty$ and $\mathbb{R}_+\times \mathcal{P}_\infty\times \mathbb{N}$ with intensity $\ddr t\otimes \mu_{\mathrm{Coag}}(\ddr \pi)$ and $\ddr t\otimes \mu_{\mathrm{Frag}}(\ddr \pi)\otimes \#(\ddr k)$ respectively. Let $\pi$ be an exchangeable random partition independent of $\mathrm{PPP}_C$ and $\mathrm{PPP}_F$. For any $n\geq 1$, set $\Pi^{n}(0)=\pi_{|[n]}$ and construct the process $(\Pi^{n}(t),t\geq 0)$ as follows:
\begin{itemize}
\item Coalescence: at an atom $(t,\pi^c)$ of $\text{PPP}_C$ such that $\pi^{c}_{|[n]}\neq 0_{[n]}$:
 \[\Pi^{n}(t)=\mathrm{Coag}(\Pi^{n}(t-),\pi^c_{|[n]}).\]
\item Fragmentation: at an atom $(t,\pi^f,k)$ of $\text{PPP}_{F}$, such that  $\pi^{f}_{|[n]}\neq 1_{[n]}$ and $k\leq n-1$, 
 \[\Pi^{n}(t)=\mathrm{Frag}(\Pi^{n}(t-),\pi^f_{|[n]},k).\]
\end{itemize}
The sequence of Markov chains $(\Pi^{n}(t),t\geq 0, n\geq 1)$ is compatible in the sense that for any $t\geq 0$ and any $m\geq n$, \[(\Pi^{m}(t)_{|[n]},t\geq 0)=(\Pi^{n}(t),t\geq 0) \text{ almost surely.}\]
This compatibility property entails the existence of a process $(\Pi(t),t\geq 0)$, taking values in the uncountable state space $\mathcal{P}_\infty$, such that  almost surely for any $n\geq 1$,
\begin{center}
$\Pi_{|[n]}(t)=\Pi^{n}(t)$ for all $t\geq 0$. 
\end{center}
Among other results, Berestycki \cite[Corollary 6, Theorem 8]{Berestycki04} has established that the process $(\Pi(t),t\geq 0)$ is a c\`adl\`ag Feller process satisfying Definition \ref{DefEFC}. 

The jump rates of the EFC process $(\Pi(t),t\geq 0)$ are prescribed by those of its restrictions $(\Pi_{|[n]},n\geq 1)$. They are easily derived from the Poisson construction in terms of $\mu_{\mathrm{Coag}}$ and $\mu_{\mathrm{Frag}}$ as follows. Let $n\in \mathbb{N}$ and $\pi\in \mathcal{P}_n$. Let $\pi^{c}$, $\pi^{f}$ be such that $\pi^{c}_{|[n]}\neq 0_{[n]}$ and $\pi^{f}_{|[n]}\neq 1_{[n]}$. Let $k\leq \#\pi$. 
\begin{itemize}
\item If $\pi\neq 1_{[n]}$, the process  $\Pi_{|[n]}$ jumps from 
$\pi$ to $\mathrm{Coag}(\pi,\pi^c)$ 
at rate: \begin{center} $\mu_{\mathrm{Coag}}(\{\pi'\in \mathcal{P}_\infty; \pi'_{|[n]}=\pi^{c}_{|[n]}\})$.\end{center}
\item If $\pi\neq 0_{[n]}$, the process  $\Pi_{|[n]}$ jumps from 
$\pi$ to $\mathrm{Frag}(\pi,\pi^f,k)$ 
at rate: \begin{center} $\mu_{\mathrm{Frag}}(\{\pi'\in \mathcal{P}_\infty; \pi'_{|[n]}=\pi^{f}_{|[n]}\})$.\end{center}
\end{itemize} 
Note that the jump rates above do no depend on the partition $\pi$.\\

The main objective of this article is to study the block-counting process $(\#\Pi(t),t\geq 0)$ and the possibility for the process to leave  the boundary $\infty$. Two behaviors  at $\infty$ are possible.

\begin{Definition}\label{dichot} Assume $\#\Pi(0)=\infty$ a.s. We say that
\begin{itemize}
\item the process stays infinite if
\begin{center}
$\forall t\geq 0; \ \#\Pi(t)=\infty$ almost surely,
\end{center} 
\item the process comes down from infinity  if 
\begin{center}
$\exists t>0; \ \#\Pi(t)<\infty$ almost surely.
\end{center}
\end{itemize}
\end{Definition}

Similarly as for pure coalescent processes, the following zero-one law holds.
\begin{lemma}[Zero-one law]\label{zeroone} Assume $\#\Pi(0)=\infty$. Set $\tau_{\infty}:=\inf\{t>0; \#\Pi(t)<\infty\}$. If $\mu_{\mathrm{Coag}}(\{\pi, \#\pi<\infty\})=0$, then either $\mathbb{P}(\tau_\infty=0)=1$ or $\mathbb{P}(\tau_\infty=\infty)=1$.    
\end{lemma}
\begin{remark}
The assumption $\mu_{\mathrm{Coag}}(\{\pi, \#\pi<\infty\})=0$ ensures that there are no coagulation events merging infinitely many blocks into finitely many. 
\end{remark}
\begin{proof}
The proof is similar as that given by Schweinsberg for pure exchangeable coalescents, see \cite[Lemma 31, p39-40]{Schweinsberg0}. We provide some details as Lemma \ref{zeroone} will play a crucial role later. The random time $\tau_\infty$ is a stopping time for the completed natural filtration of $(\Pi(t),t\geq 0)$. Since $(\Pi(t),t\geq 0)$ is a Feller process, then by Blumenthal's zero-one law, one has $\mathbb{P}(\tau_{\infty}=0)\in \{0,1\}$. It remains to show that the event $\{\tau_\infty \in (0,\infty)\}$ has probability zero. Consider first the event
\[
\{0<\tau_\infty<\infty,\ \#\Pi(\tau_\infty-)=\infty \text{ and }\#\Pi(\tau_\infty)<\infty\}.\]
On this event, $\tau_\infty$ must be an atom of $\text{PPP}_{C}$ at which infinitely many blocks merge into finitely many. Since by assumption $\mu_{\mathrm{Coag}}(\{ \pi\in \mathcal{P}_\infty, \#\pi<\infty\})=0$ and all partitions atoms of $\mathrm{PPP}_C$ have infinitely many blocks, the latter event has probability zero. We now show that the event \[\{\tau_\infty\in (0,\infty), \ \#\Pi(\tau_\infty-)<\infty\}\] has also probability zero. 
For any $b\in \mathbb{N}$, set $\lambda_{b+1}:=\mu_{\mathrm{Coag}}(\{\pi; \pi_{|[b+1]}\neq 0_{[b+1]}\})<\infty$, the rate at which a coalescence involving $b+1$ blocks occurs. Fix $b$ and $n_1<n_2<\cdots<n_b$ in $\mathbb{N}$. Consider the event
\[I(b,n_1,\cdots,n_b):=\{\tau_\infty\in (0,\infty), \ \#\Pi(\tau_\infty-)=b, \min \Pi_i(\tau_\infty-)=n_i, \text{for all } i\in [b]\}.\]
Let $T_0=0$, and choose an integer $p_1\geq 2$ such that $p_1 \underset{\Pi(0)}{\not \sim} n_i$ for all $i\in [b]$ (such a $p_1$ exists since $\#\Pi(0)=\infty$). Let $T_1:=\inf \{t>T_0; p_1 \underset{\Pi(t)}{\sim} n_i \text{ for some }i \in [b]\}$, this is a coalescence time between the particular block containing $p_1$ and at most $b$ other possible blocks. The rate of $T_1$ is thus at most $\lambda_{b+1}$. Since by assumption, infinitely many blocks cannot coagulate into finitely many by a single jump, $T_1<\tau_\infty$ a.s and thus $\#\Pi(T_1)=\infty$ a.s. Recursively, we can choose an integer $p_m$ such that $p_m \underset{\Pi(T_{m-1})}{\not \sim} n_i$ for all $i\in [b]$ and define \[T_{m}:=\inf \{t>T_{m-1}; p_m \underset{\Pi(t)}{\sim} n_i \text{ for some }i \in [b]\}.\] By construction $T_{m-1}\leq T_{m}< \tau_\infty$ a.s. for any $m\geq 1$ and $T_{m}-T_{m-1}$ is stochastically greater than an exponential random variable with parameter $\lambda_{b+1}$. We deduce that on $I(b,n_1,\cdots, n_b)$, $\tau_\infty\geq \sum_{m=1}^{\infty}(T_{m}-T_{m-1})= \underset{m\rightarrow \infty}{\lim} T_{m}=\infty \text{ a.s.}$ Therefore,  $I(b,n_1,\cdots, n_b)$ has probability zero, as well as the event  
\[\{\tau_\infty\in (0,\infty), \ \#\Pi(\tau_\infty-)<\infty\}=\bigcup_{b\in \mathbb{N}, n_{1}<n_2<\cdots<n_b}I(b,n_1,\cdots,n_b).\] 
We conclude that $\mathbb{P}(\tau_\infty \in (0,\infty))=0$.
\end{proof}
\begin{remark}\label{card}
We highlight that the map $\#: \mathcal{P}_\infty \rightarrow \mathbb{N}\cup \{\infty\}$ is not continuous with respect to $d$. For instance, if for any $k\in \mathbb{N}$, $\pi^{(k)}:=\{[k],\{k+1\},\cdots\}$, then $\underset{k\rightarrow \infty}{\lim} d(\pi^{(k)},1_{\mathbb{N}})=\underset{k\rightarrow \infty}{\lim} 1/k=0$ and $\pi^{(k)}$ converges to $1_{\mathbb{N}}$, even though $\#\pi^{(k)}=\infty$ and $\#1_{\mathbb{N}}=1$. However, it is important to note that if $\#\pi=\infty$ and $\underset{k\rightarrow \infty}{\lim} d(\pi^{(k)},\pi)=0$, then $(\#\pi^{(k)},k\geq 1)$ converges towards $\#\pi=\infty$, as $k$ goes to $\infty$. Indeed, recalling \eqref{distance}, if $\underset{k\rightarrow \infty}{\lim} d(\pi^{(k)},\pi)=0$, then $n_k:=\max\{n\geq 1, \pi^{(k)}_{|[n]}=\pi_{|[n]}\} \underset{k\rightarrow \infty}{\longrightarrow} \infty$ and the result follows since for any $k$, $\#\pi^{(k)}\geq \#\pi_{|[n_k]}$.
\end{remark}

The remark above indicates that it will be necessary to shed some light on the Markov property of the process  $(\#\Pi(t),t\geq 0)$, as well as on the regularity of its paths.  We postpone this discussion for  simple EFC processes to Section \ref{simpleEFC}, see the forthcoming Proposition \ref{generatorlemma} and Remark \ref{remarkonmarkov2}.

\subsection{Exchangeable coalescent processes and their number of blocks}\label{cdicoal}
Pure exchangeable coalescent processes, namely EFC processes with $\mu_{\mathrm{Frag}}\equiv 0$, have received a lot of attention. We refer to the seminal papers of Pitman \cite{Pitman99}, Schweinsberg \cite{Schweinsberg0}, Sagitov \cite{MR1742154} and M\"ohle and Sagitov \cite{MR1880231}. See also Berestycki's book \cite[Chapters 3 and 4]{Beresbook} for a recent account on fine properties of the so-called $\Lambda$-coalescents.

It is worth noticing that the number of blocks in any pure exchangeable coalescent has decreasing sample paths. 
If the coalescent comes down from infinity, in the sense of Definition \ref{dichot}, then it stays finite a.s. after it has comed down. This is a striking difference with the block counting process of an EFC process whose sample paths are not monotone.

Only sufficient conditions entailing coming down from infinity are known for a general coagulation measure $\nu_{\mathrm{Coag}}$ carried over $\mathcal{P}_{\mathrm{m}}$ ($\Xi$-coalescents), see Herriger and M\"ohle \cite{MohleHerriger}. For the sake of simplicity, we focus now on coalescences in which there are no simultaneous multiple collisions and no coagulation of all blocks at once. Namely those satisfying \begin{equation}\label{lambdacoal}
c_{\mathrm{k}}\geq 0,\ \nu_{\mathrm{Coag}}(\{s\in\mathcal{P}_{\mathrm{m}}; \ s_2>0\})=0, \ \nu_{\mathrm{Coag}}(\{s\in\mathcal{P}_{\mathrm{m}}; \ s_1=1\})=0.
\end{equation} 
Since $\nu_{\mathrm{Coag}}$ is carried over $\{s\in \mathcal{P}_{\mathrm{m}}; s_2=0\}$, the measure $\nu_{\mathrm{Coag}}$ can be considered as a measure on $[0,1]$ and the atoms of the Poisson point process $\text{PPP}_C$ have partitions with only one non-singleton block.  When moreover $c_{\mathrm{k}}=0$, it is often useful to describe the coalescent part of Section \ref{EFCdef} as follows.  Associate to each atom $(t,\pi^{c})$ of $\text{PPP}_C$, the sequence of random variables $(X_k)_{k\geq 1}$ defined by
\[X_k=1 \text{ if } \{k\} \text{ is not a singleton block of } \pi^{c} \text{ and } X_k=0 \text{ otherwise}.\] 
By assumption \eqref{lambdacoal} and by definition, we have \[k\underset{\pi^{c}}{\sim} \ell  \Longleftrightarrow X_k=X_{\ell}=1.\]
Given $|\pi^c|^{\downarrow}=x\in (0,1)$, $(X_k)_{k\geq 1}$ is a sequence of i.i.d Bernoulli random variables with parameter $x$.  The coalescence event occurring at time $t$:
\[\Pi_{|[n]}(t)=\mathrm{Coag}(\Pi_{|[n]}(t-),\pi^c_{|[n]})\]
can now be described as follows:
all blocks of $\Pi_{|[n]}(t-)$ whose index $k\leq \#\Pi_{|[n]}(t-)$ satisfies $X_k=1$,  merge together.  We refer for instance to \cite[Theorem 3.2 and Corollary 3.1]{Beresbook}, in particular to see how to incorporate binary coalescences when $c_{\mathrm{k}}>0$. \\

A coalescent process whose coagulation measure satisfied \eqref{lambdacoal} is often called in the literature a $\Lambda$-coalescent. The prefix $\Lambda$ stands for the finite measure $\Lambda:=c_{\mathrm{k}}\delta_0+x^{2}\nu_{\mathrm{Coag}}(\ddr x)$, which characterizes the law of the process. More precisely, the process $(\Pi(t),t\geq 0)$ is characterized in law by the jump rates of its restrictions, namely by the sequence $(\lambda_{n,k}, 2\leq k\leq n)_{n\geq 2}$ defined by
\begin{align}\label{lambdank}
\lambda_{n,k}&:=\mu_{\mathrm{Coag}}(\{\pi; \text{ the non-singleton block of } \pi_{|[n]} \text{ has } k \text{ elements}\})\nonumber \\
&=c_{\mathrm{k}}\mathbbm{1}_{\{k=2\}}+\int_{0}^{1}x^{k}(1-x)^{n-k}\nu_{\mathrm{Coag}}(\ddr x).
\end{align}
As recalled in the Introduction, Schweinsberg \cite{Schweinsberg00} has established a necessary and sufficient condition for coming down from infinity of $\Lambda$-coalescents. Recall $\Phi(n)$ defined in \eqref{phi}. Some binomial calculations, see \cite{Schweinsberg00}, yield the following other expression of $\Phi(n)$.  For any $n\geq 2$,  
\begin{equation} \label{phi2} \Phi(n)=c_{\mathrm{k}}\binom{n}{2}+\int_{0}^{1}\left(nx-1+(1-x)^{n}\right)\nu_{\mathrm{Coag}}(\ddr x). \end{equation} 
It is not difficult to verify, from this identity, that $(\Phi(n)/n ,n\geq 1)$ is non-decreasing. One can also check analytically that
$\Phi(n)\underset{n\rightarrow \infty}{\sim} \Psi(n)$
with $\Psi$ the function of the L\'evy-Khintchine form :
\begin{equation} \label{psi} 
\Psi(u)=\frac{c_{\mathrm{k}}}{2}u^{2}+\int_{0}^{1}\left(e^{-xu}-1+ux\right)\nu_{\mathrm{Coag}}(\ddr x).
\end{equation}
See for instance \cite[Lemma 3.3]{MohleHerriger}. A first consequence of this latter equivalence is that the condition of coming down from infinity \eqref{CDI} is equivalent to the integrability condition $\int_2^{\infty}\frac{\ddr u}{\Psi(u)}<\infty$. See \cite[Section 4.3]{Beresbook} for some probabilistic interpretation of this equivalence in the stable case.

One other interest in using $\Psi$ instead of $\Phi$ is that we can easily apply Tauberian theorems to find an equivalent of $\Psi$ (and then of $\Phi$) when the measure $\nu_{\mathrm{Coag}}$ has some properties of regular variation. The following is a direct application of the Tauberian theorem (see e.g. \cite[Page 10]{Ber96}). For any $x\in (0,1]$, set $\bar{\nu}_{\mathrm{Coag}}(x):=\nu_{\mathrm{Coag}}([x,1])$ and $\bar{\bar{\nu}}_{\mathrm{Coag}}(x):=\int_{x}^{1}\bar{\nu}_{\mathrm{Coag}}([v,1])\ddr v$. 
If $\bar{\bar{\nu}}_{\mathrm{Coag}}(x)\underset{x\rightarrow 0}{\sim}x^{\rho-1}\frac{L(x)}{\Gamma(\rho)}$ for some $\rho \in (0,\infty)$ and $L$ is a slowly varying function then 
\begin{equation}\label{equivpsi}
\Psi(u)\underset{u\rightarrow \infty}{\sim} u^{2-\rho}L(1/u).
\end{equation}
For instance, assume that $\nu_{\mathrm{Coag}}(\ddr x)=f(x)\ddr x$ with $f$ such that $f(x)x^{2+\beta}\underset{x\rightarrow 0}{\longrightarrow} c>0$ with $\beta\in (0,1)$, then $\bar{\bar{\nu}}_{\mathrm{Coag}}(x)\underset{x\rightarrow 0}{\sim} \frac{c}{\beta(\beta+1)}x^{-\beta}$ and by taking $\rho=1-\beta$ and $L(x)=\Gamma(1-\beta)c$ for all $x\in [0,1]$ in \eqref{equivpsi}, one gets $\Psi(n)\underset{n\rightarrow \infty}{\sim} c\frac{\Gamma(1-\beta)}{\beta(\beta+1)} n^{\beta+1}$ and therefore
\[\Phi(n)\underset{n\rightarrow \infty}{\sim} dn^{\beta+1}\] 
with  $d=c\frac{\Gamma(1-\beta)}{\beta(\beta+1)}$. 
Applying now the Tauberian theorem with $\rho=1$ and $L(x)=\log(1/x)^{\beta}$ gives that any coagulation measure $\nu_{\mathrm{Coag}}$ for which $\bar{\bar{\nu}}_{\mathrm{Coag}}(x)\underset{x\rightarrow 0}{\sim} c\log(1/x)^{\beta}$ satisfies $$\Phi(n)\underset{n\rightarrow \infty}{\sim} cn(\log n )^{\beta}.$$
The conditions that bear on the function $\Phi$ of Proposition \ref{regularcdi} and Proposition \ref{slowcdi} are therefore satisfied for the coagulation measures constructed above.
\begin{remark}\label{crazyexample}  In general, the L\'evy-Khintchine function $\Psi$ may have different upper and lower indices at $\infty$. We refer to Bertoin's lecture notes \cite[Chapter 5]{MR1746300} for the definition of these indices and for seeing how to construct a L\'evy measure $\nu_{\mathrm{Coag}}$ providing  such a function $\Psi$. In this case, the parameters $\theta^{\star}$ and $\theta_\star$ might not coincide, but 
we shall not consider further this question here.
\end{remark}

\subsection{Simple EFC processes}\label{simpleEFC}
Recall \eqref{lambdacoal} and its meaning in terms of coalescence. We now introduce the so-called simple EFC processes. Since any block in an exchangeable random partition of $\mathbb{N}$ is either singleton or infinite, the condition \eqref{infiniteblock} is equivalent to the assumption that $\mu_{\mathrm{Frag}}$  is supported by partitions with no singletons. 

\begin{Definition} An EFC process is called simple if its coagulation measure satisfies \eqref{lambdacoal} and if its fragmentation measure has finite total mass and is supported by partitions with no singletons. 
\end{Definition}
\begin{remark} The name coined ``simple" follows Bertoin's terminology for the $\Lambda$-coalescents, see \cite[Section 4.4]{coursbertoin}. Beside the fact that there  is no formation of dust, a simple EFC process has no simultaneous multiple coagulations and can only fragmentate a single block at a time.
\end{remark}
According to \eqref{represmufrag}, the first assumption on the fragmentation measure, $\mu_{\mathrm{Frag}}(\mathcal{P}_\infty)<\infty$, is equivalent to 
$c_e=0 \text{ (no erosion coefficient) and }  \nu_{\mathrm{Disl}}(\mathcal{P}_{\mathrm{m}})<\infty.$
The second assumption on its support \eqref{infiniteblock} is equivalent to having a dislocation measure supported by $\cup_{k\in \mathbb{N}\cup \{\infty\}}\mathcal{P}^{k}_{\mathrm{m}}$ where for any $k\in \mathbb{N}\cup \{\infty\}$
 \[\mathcal{P}_{\mathrm{m}}^{k}:=\left\{s\in \mathcal{P}_{\mathrm{m}} ;s_i>0,\ \forall i\in [k+1],\ \sum_{i=1}^{k+1}s_i=1\right\}.\] 
By Kingman's paint-box representation, see Figure \ref{paintbox}, if $\pi^{f}$ is an atom of $\text{PPP}_F$ then, on the event $\{|\pi^{f}|^{\downarrow}\in \mathcal{P}^{k}_{\mathrm{m}}\}$, one has $\#\pi^{f}=k+1$ almost surely. We stress that we allow the value $k=\infty$, so that fragmentation into infinitely many pieces are possible.


Since $\mu_{\mathrm{Frag}}$ is assumed to be finite, the process $(\#\Pi(t),t\geq 0)$ restricted to $\mathbb{N}$ evolves as a classical continuous-time process with no instantaneous integer state. In order to take into account instantaneous coming down from infinity and possible explosion, we consider as state-space, the one-point compactification of $\mathbb{N}$, which we denote by $\bar{\mathbb{N}}$. The restricted state-space $\mathbb{N}$ is thus endowed with the discrete topology, and for any $m\geq 1$, $\{\infty\}\cup [|1,m|]^{c}$  forms a neighborhood of $\infty$. 

\begin{proposition}\label{generatorlemma} Let $(\Pi(t),t\geq 0)$ be a simple EFC process. The block-counting process $(\#\Pi(t),t\geq 0)$ is a right-continuous process valued in $\bar{\mathbb{N}}$.  Moreover, at any time $t>0$ such that $\#\Pi(t-)<\infty$, $\underset{h\rightarrow 0^{+}}{\lim} \#\Pi(t-h)=\#\Pi(t-)$ a.s. 
The process $(\#\Pi(t),t<\zeta)$ started from $n$ and  stopped at its first explosion time $\zeta:=\inf\{t>0; \#\Pi(t-)=\infty\ \mathrm{or}\ \#\Pi(t)=\infty\}$ is Markov and its generator $\mathcal{L}$ acts on 
\begin{equation}\label{domain} \mathcal{D}:=\left\{g:\bar{\mathbb{N}}\rightarrow \mathbb{R};\ \forall n\in \bar{\mathbb{N}},\sum_{k\in \mathbb{N}\cup \{\infty\}}|g(n+k)|\mu(k)<\infty\right\},
\end{equation}
as follows: for any $g\in \mathcal{D}$ 
 \[\mathcal{L}g:=\mathcal{L}^{c}g+\mathcal{L}^{f}g\]
with, for any $n\in \mathbb{N}$
\begin{equation}\label{coalpart}
\mathcal{L}^{c}g(n):=\sum_{k=2}^{n}\binom{n}{k}\lambda_{n,k}[g(n-k+1)-g(n)]
\end{equation}
and
\begin{equation}\label{fragmentpart}
\mathcal{L}^{f}g(n):=n\sum_{k=1}^{\infty}\mu(k)[g(n+k)-g(n)]+n\mu(\infty)[g(\infty)-g(n)]. 
\end{equation} 
\end{proposition}
\begin{remark}
Since the splitting measure $\mu$ is finite, $\mathcal{D}$ contains all bounded functions on $\bar{\mathbb{N}}$. In particular, any continuous function on $\bar{\mathbb{N}}$, i.e. any function $g$ defined on $\mathbb{N}\cup\{\infty\}$ such that $\underset{n\rightarrow \infty}{\lim} g(n)=g(\infty)<\infty$ belongs to $\mathcal{D}$. Moreover, if $\mu(\infty)=0$ then $\mathcal{D}$ contains all function $g$ defined over $\mathbb{N}$ such that the series $\sum_{k=1}^{\infty}|g(n+k)|\mu(k)$ converges for all $n\in \mathbb{N}$. 
\end{remark}
\begin{proof} 
By assumption, $\mu_{\mathrm{Frag}}(\mathcal{P}_\infty)<\infty$, and the fragmentations occur at finite rate. Moreover starting from a partition $\pi$ with finitely many blocks, there are only finitely many possibilities of coagulations, each with a finite rate. 
Partitions with finitely many blocks are therefore holding points for the process $(\Pi(t),t\geq 0)$. We now study the left-limits. Let $t>0$ and  assume that $\#\Pi(t-)<\infty$. By the assumption \eqref{lambdacoal}, blocks cannot all coagulate at once. Since $\#\Pi(t-)<\infty$, the partition reached by the process before the state $\Pi(t-)$ has also finitely many blocks. Hence, almost surely, for $h$ small enough $\Pi(t-h)=\Pi(t-)$ and therefore, $\underset{h\rightarrow 0}{\lim} \#\Pi(t-h)=\#\Pi(t-)$. 
The right-continuity at any time $t$ such that $\#\Pi(t)<\infty$ is obtained along similar arguments, using the fact that the process stays an exponential time at the partition $\Pi(t)$. Assume now $\#\Pi(t)=\infty$, since $\Pi$ is right-continuous, $d\big(\Pi(t+h),\Pi(t)\big)\underset{h\rightarrow 0}{\longrightarrow} 0$, and one has $\#\Pi(t+h)\underset{h\rightarrow 0}{\longrightarrow}  \#\Pi(t)=\infty$, see Remark \ref{card}.

The form of the generator will be deduced from the Poisson construction. The part with negative jumps corresponds to the generator of the number of blocks in a $\Lambda$-coalescent started from a partition with $n$ blocks. We refer for instance to \cite[page 203]{coursbertoin} and focus on the positive jumps. Let $n\in \mathbb{N}$ and assume $\#\Pi(0)=n$. Let $m\in \bar{\mathbb{N}}$. Consider an atom $(t,\pi^{f},j)$ of $\text{PPP}_F$. The atom $(t,\pi^{f},j)$ is seen by $\Pi_{|[m]}(t-)$ if $\pi_{|[m]}^{f}\neq 1_{[m]}$ and $j\leq \#\Pi_{|[m]}(t-)$. By definition of the fragmentation operator, see Definition \ref{coagfrag}, we have
\begin{equation}\label{cardfrag} \#\Pi_{|[m]}(t)=\#\Pi_{|[m]}(t-)-1+\#\{\Pi_{j}(t-)\cap \pi^{f}_i\cap [m],1\leq i\leq \#\pi^{f}\}.\end{equation}
Let $m=\infty$, first observe that if $\#\Pi(t-)=\infty$ then $\#\Pi(t)=\infty$ and $t$ is not a jump time (by independence of the Poisson point processes it cannot be a coalescence time).  Note that when $\#\Pi(t-)<\infty$, each block of the partition $\Pi(t-)$ is infinite a.s, otherwise the partition $\Pi(t)$, just after the fragmentation, would contain blocks of finite size, which is not possible, since $\Pi(t)$ is exchangeable and has no dust. 

Notice that by the Poisson construction, $\Pi_j(t-)$ is independent of $\pi^{f}$.  Assume first $\#\pi^{f}=k+1$ with $k\in \mathbb{N}$. So that, for any $i\in [k+1]$, $|\pi^{f}_i|>0$  almost surely. Let $(U_{\ell}, \ell\geq 1)$ be a sequence of i.i.d uniform random variables on $[0,1]$, independent of $\Pi_j(t-)$. By the paintbox representation of $\pi^f$, and since $\#\Pi_j(t-)=\infty$, one has for any $1\leq i\leq \#\pi^f$
\begin{equation*}\label{itsbranching} \mathbb{P}\big(\Pi_{j}(t-)\cap \pi_i^{f}=\emptyset\big)=\mathbb{P}\big(\forall \ell \in \Pi_j(t-), U_{\ell} \text{ is not in a subinterval of length } |\pi^{f}_i|\big)=0. 
\end{equation*}
Therefore $\#\{\Pi_{j}(t-)\cap \pi^{f}_i, 1\leq i\leq \#\pi^{f}\}=\#\pi^{f}=k+1$, and the right hand side of \eqref{cardfrag} equals $n+k$ a.s. Recall the intensity of $\text{PPP}_F$, $\mu_{\mathrm{Frag}}(\ddr \pi)\otimes \#(\ddr j)$. When $\#\Pi(t-)=n$, the process jumps to $n+k$ at rate $$n\mu_{\mathrm{Frag}}(\{\pi, \ \#\pi^{f}=k+1\})=n\nu_{\mathrm{Disl}}(\mathcal{P}^{k}_{\mathrm{m}})=n\mu(k).$$ 
Suppose now $\#\pi^{f}=\infty$. By assumption \ref{simpleEFC},  $|\pi^{f}_i|>0$ for all $i\geq 1$, and the same argument shows that  the r.h.s in \eqref{cardfrag} is infinite. In this case, the process jumps from $n$ to $\infty$. Finally, one sees that the rate of jump from $n$ to $\infty$ is given by $n\nu_{\mathrm{Disl}}(\mathcal{P}^{\infty}_{\mathrm{m}})$ where we recall $\mathcal{P}^{\infty}_{\mathrm{m}}:=\{\mathrm{s}\in \mathcal{P}_{\mathrm{m}};\ s_i>0 \text{ for all } i\geq 1 \text{ and }  \sum_{i=1}^{\infty}s_i=1\}$. 
\end{proof}
\begin{remark}\label{remarkonmarkov1}
When $m<\infty$, the process $(\#\Pi_{|[m]}(t),t\geq 0)$ is not Markov in general. Indeed, it may occur that for some index $j$, $\Pi_{j}(t-)\cap \pi^{f}_i\cap [m]=\emptyset$ for some $i$ and therefore that $\#\Pi_{|[m]}(t)$ depends on the constituent elements of the blocks  of $\Pi_{|[m]}(t-)$ and not only on its number of blocks.
\end{remark}
\begin{remark}\label{remarkonmarkov2} Let $t>0$. Conditionally on $\#\Pi(t)<\infty$, if one denote by $\zeta\circ \theta_t$ the first explosion time (possibly infinite) of the process $(\#\Pi(t+s), s>0)$, then Proposition \ref{generatorlemma} entails that $(\#\Pi(t+s),  s<\zeta\circ \theta_t)$ is Markov and has the same law as $(\#\Pi(s),s<\zeta)$ when started from $n:=\#\Pi(t)<\infty$. The Markov property of $(\#\Pi(t),t\geq 0)$ is however more involved at the times $t$ such that $\#\Pi(t)=\infty$. Indeed, if by the Markov property of the EFC process, for any such time $t$,  conditionally on $\mathcal{F}_t$ (the natural filtration), $(\Pi(s+t),s\geq 0)$  has a law only depending on $\Pi(t)$, a priori the law of the process $(\#\Pi(t+s),s\geq 0)$ could depend on the shape of the blocks of $\Pi(t)$. We shall circumvent this difficulty by only making use of the Markov properties of the processes $(\Pi(t),t\geq 0)$ and $(\#\Pi(t),t<\zeta)$. We refer the reader interested in this question to Theorem 3.7 in Foucart and Zhou \cite{FoucartZhouWF}, where it is established that $(\#\Pi(t),t\geq 0)$ is actually a Feller Markov process on $\bar{\mathbb{N}}$, the one-point compactification of $\mathbb{N}$.\end{remark}

\section{Proof of Theorem \ref{mainthm}}\label{proof}
We start by presenting a heuristic argument enlightening Theorem \ref{mainthm}. Recall that we work with a coalescence measure $\Lambda$ satisfying $\Lambda(\{1\})=0$ and Schweinsberg's condition \eqref{CDI}: $\sum_{k=2}^{\infty}\frac{1}{\Phi(k)}<\infty$. This entails that the pure $\Lambda$-coalescent comes down from infinity instantaneously.
The proof of Theorem \ref{mainthm} is mainly based on martingale arguments involving the bounded function \[f: n\mapsto \sum_{j=n+1}^{\infty}\frac{1}{\Phi(j)}.\] This function appears in several previous works. We refer the reader to Bertoin's book \cite[Proposition 4.9, page 202]{coursbertoin} and Berestycki's book \cite[Chapter 3, page 70]{Beresbook}. \\

\textbf{Heuristic}. Recall the definition of $\Phi$ in \eqref{phi}, its meaning and the definitions of $\theta^{\star}$ and $\theta_{\star}$ in \eqref{theta}.  Although typically the process $(\#\Pi(t),t\geq 0)$ does not only move down one level at a time  (except in the pure Kingman case), $f(n)$ turns out to be a rather sharp measure of the time needed for the pure coalescent process to go from infinity to level $n$ when  $n$ is large. Indeed, it is actually established in \cite{Beres10} that the speed of coming down from infinity of the $\Lambda$-coalescent is precisely the inverse function of $f$. Keeping this in mind, the  conditions $\theta^{\star}<1, \theta_{\star}>1$ for coming down from infinity and staying infinite, can be understood as follows: let $Z$ be the number of new blocks after a fragmentation, this is  a random variable with law $\mu(\cdot)/\mu(\bar{\mathbb{N}})$. A simple application of Fubini's theorem  ensures that for any $n\in \mathbb{N}$, \[\sum_{k=1}^{\infty}\frac{n\bar{\mu}(k)}{\Phi(n+k)}=n\mu(\bar{\mathbb{N}}) \mathbb{E}\left[\sum_{k=n+1}^{n+Z}\frac{1}{\Phi(k)}\right]=\frac{\mathbb{E}[f(n)-f(n+Z)]}{1/\mu(\bar{\mathbb{N}})n}.\] Therefore, the  condition $\theta^{\star}<1$ (respectively $\theta_{\star}>1$) can be seen  as the condition under which for all large enough $n$, the mean time for the $\Lambda$-coalescent to go from $n+Z$ to $n$, is stricly smaller (respectively larger) than $1/\mu(\bar{\mathbb{N}})n$, the mean arrival time of a typical fragmentation (whose size is  $Z$) when there are $n$ blocks.\\

The proof of Theorem \ref{mainthm} will not make use of the speed of coming down from infinity previously mentioned. We establish the first part of Theorem \ref{mainthm}, namely that $(\Pi(t), t\geq 0) $ comes down from infinity when $\theta^{\star}<1$, in Subsection \ref{CDIsubsec}. The second part of Theorem \ref{mainthm} is shown in Subsection \ref{Stayinfinitesubsec}. We start by preliminaries on the process $(\#\Pi(t),t\geq 0)$.
\subsection{Preliminaries}
We first explain why we can restrict our study to simple EFC processes starting from partitions with blocks of infinite size. 
\begin{lemma}\label{nodust} Let $(\Pi(t),t\geq 0)$ be a simple EFC process whose coalescence measure $\Lambda$ satisfies $\int _0^1 x^{-1}\Lambda(\ddr x)=\infty$. Then, the process $(\Pi(t),t\geq 0)$ started from any exchangeable partition $\Pi(0)$, has no dust  i.e. no singleton blocks, at any strictly positive time, almost surely. \end{lemma}
\begin{proof}  The assumption $\int_0^1 x^{-1}\Lambda(\ddr x)=\infty$ ensures that the pure $\Lambda$-coalescent process $(\Pi^{C}(t),t\geq 0)$ has no dust. Namely at any time $t>0$, $\Pi^{C}(t)$ has no singleton blocks, even though $\Pi^{C}(0)$ has some. We refer to \cite[Theorem 3.5]{Beresbook}. Note also that since the fragmentation measure has finite mass and is carried over partitions with no singleton blocks, the pure fragmentation process $(\Pi^{F}(t),t\geq 0)$, started from a partition with no singleton blocks,  has no singleton blocks at any time $t\geq 0$. 

We now follow the arguments of the proof of \cite[Proposition 16]{Berestycki04}, see Section 5.2, page 794-795. Consider the partition-valued process $(\Pi(t),t\geq 0)$ constructed from $\mathrm{PPP}_C$ and $\mathrm{PPP}_F$ as follows: at any fragmentation time $(t,\pi^f,k)$, the partition $\Pi(t)$ is obtained from $\Pi(t-)$ as for a classical EFC process; and at any coalescence time $(t,\pi^c)$, the partition $\Pi(t)$ is obtained by merging the blocks $\Pi_i(t-)$ and $\Pi_j(t-)$ if and only if $\min \Pi_i(t)$ and $\min \Pi_j(t-)$ (instead of $i$ and $j$) belong to the same block of $\pi^c$. In other words, recalling the definition of the i.i.d sequence $(X_k,k\geq 1)$, see Section \ref{cdicoal}, at any coalescence time, the block with index $k$ takes part to the coalescence if and only if
\[\tilde{X}_k:=X_{\min \Pi_k(t-)}=1.\]   
By independence of $\Pi_{|[n]}(t-)$ and $\pi^c_{|[n]}$, and exchangeability of $\pi^c_{|[n]}$,  for any $n\geq 1$, we see that $(\tilde{X}_k, 1\leq k\leq n)$ has the same law as $(X_k, 1\leq k\leq n)$. The process $(\Pi(t),t\geq 0)$ is thus a simple EFC process with fragmentation measure $\mu_{\mathrm{Frag}}$ and coagulation measure $\Lambda$. 
 
Denote by $(\Pi^{C}(t),t\geq 0)$ the pure coalescent process obtained  from $\mathrm{PPP}_C$ by the modified Poisson construction as above. By definition the singleton block $\{i\}$ takes part to a coalescence event at time $t$ in $\Pi^c$ if and only if $\tilde{X}_{j}=1$ where $j$ is such that $\Pi^{C}_j(t-)=\{i\}$. By definition $\tilde{X}_{j}=X_i$, thus $X_i=1$. Similarly the singleton block $\{i\}$ takes part to a coalescence event at time $t$ in the EFC process $\Pi$ if and only if $\tilde{X}_{j}=1$ where $j$ is such that $\{i\}=\Pi_j(t-)$. Hence $\tilde{X}_{j}=X_i$ and $X_i=1$. The coalescences of singletons are therefore coupled for both processes $(\Pi(t),t\geq 0)$ and $(\Pi^{C}(t),t\geq 0)$. Since there is no singletons in $\Pi^C$, there is also none in $\Pi$. \end{proof}
\noindent Recall that we work under the assumptions \eqref{CDI} and $\Lambda(\{1\})=0$. The pure $\Lambda$-coalescent comes down from infinity instantaneously and has no singleton blocks. By Lemma \ref{nodust},  for any $t>0$, the blocks of the EFC process at time $t$, $\Pi(t)$, are also of infinite size. Moreover  Lemma \ref{zeroone} ensures that the process $\Pi$ stays infinite if and only if the EFC $(\Pi(t+s),s\geq 0)$ stays infinite for all $t>0$. We can therefore suppose without loss of generality that $\#\Pi(0)=\infty$ and that for any $i\in \mathbb{N}$, $\#\Pi_i(0)=\infty$ a.s. From now on, we work with such an initial partition.\\

We will study the  process $(\#\Pi(t),t\geq 0)$ from a monotone coupling $(\Pi^{(n)}(t),t\geq 0)$ of $(\Pi(t),t\geq 0)$ satisfying $\#\Pi^{(n)}(0)=n$ and $\#\Pi^{(n)}(t)\leq \#\Pi^{(n+1)}(t)$ for all $t$ a.s. See the forthcoming Lemma \ref{coupling1}. Note that it is not so obvious, at a first glance, that such a monotone coupling exists. Firstly observe that for any partitions $\pi,\pi'\in \mathcal{P}_\infty$ and any $k\geq 1$, the operators $\mathrm{Coag}(\cdot,\pi)$ and $\mathrm{Frag}(\cdot,\pi,k)$, see Definition \ref{coagfrag}, still make sense when acting on a collection of disjoint subsets that are indexed in the order of their least elements: if $B=(B_1,\cdots, B_n)$ , for some $n\in \bar{\mathbb{N}}$, are disjoint subsets, with for any $i\leq j$, $\min B_i \leq \min B_j$, then the following are well-defined
\[\mathrm{Coag}(B,\pi^{c}):=\left\{ \cup_{j\in \pi^c_i}B_j, i\geq 1\right\} \text{ and }\mathrm{Frag}(B,\pi^{f},k):=\left\{B_j, j\neq k, B_k\cap \pi^{f}_i, i\geq 1\right\}^{\downarrow}.\]
As for partitions, $\#B$ denotes the number of non-empty subsets of any collection $B$ and by convention when $\#B<\infty$, $B_i:=\emptyset$ for any $i>\#B$.  We also stress that the space of ordered collections of disjoint subsets of $\mathbb{N}$, call it $\mathcal{S}$, is embedded in the space of partitions $\mathcal{P}_\infty$. Indeed, if $B\in \mathcal{S}$, then  $\pi:=\left\{B_1,\cdots, B_n, \left(\cup_{i=1}^{n}B_i\right)^{c}\right\}^{\downarrow} \in \mathcal{P}_\infty$. \\

The following simple lemma will play a role in our first coupling argument.\begin{lemma}\label{cardcoal} 
Let $\pi^{c}\in \mathcal{P}_\infty$. Let $B^{1}$ and $B^{2}$ be two finite collections of disjoint subsets of $\mathbb{N}$ ordered by their least element. If $\#B^{1}\leq \#B^{2}$ then \[\#\mathrm{Coag}(B^{1},\pi^{c})\leq \#\mathrm{Coag}(B^{2},\pi^{c}).\]
\end{lemma}
\begin{proof}
Assume by contradiction that $\#\mathrm{Coag}(B^{1},\pi^{c})>\#\mathrm{Coag}(B^{2},\pi^{c})$. Let $i:=\#\mathrm{Coag}(B^{1},\pi^{c})$. On the one hand, $\cup_{j\in \pi^{c}_i}B^{2}_j=\emptyset$ and thus for all $j\in \pi^{c}_i$, $B^{2}_j=\emptyset$. Therefore $\#B^2<\min \pi_i^{c}$. On the other hand, by definition $\cup_{j\in \pi^{c}_i}B^{1}_j\neq \emptyset$ and there exists $j\in \pi_i^{c}$ such that $B^{1}_j\neq \emptyset$. This entails $\#B^{1}\geq j\geq \min \pi_i^{c}$ and leads to the contradiction $\#B^{1}>\#B^{2}$. 
\end{proof}

Recall the Poisson point processes $\mathrm{PPP}_F$ and $\mathrm{PPP}_C$. Let $n\geq 1$. We now construct a process $(\Pi^{(n)}(t),t\geq 0)$ valued in $\mathcal{S}$, started from $(\Pi_1(0),\cdots, \Pi_n(0))$, which follows all fragmentations and coagulations involving integers belonging to $\cup_{i=1}^{n}\Pi_i(0)$.\\

For any $m\in \mathbb{N}$, set $\Pi^{(n),m}(0):=(\Pi_1(0)\cap [m],\cdots, \Pi_n(0)\cap [m])$ and 

\begin{itemize}
\item if $(t,\pi^{c})$ is an atom of $\mathrm{PPP}_C$ such that $\pi^{c}_{|[m]}\neq 0_{[m]}$, then 
\[\Pi^{(n),m}(t)=\mathrm{Coag}(\Pi^{(n),m}(t-),\pi^{c}_{|[m]}),\]
\item if $(t,\pi^{f},k)$ is an atom of $\mathrm{PPP}_F$ such that $\pi^{f}_{|[m]}\neq 1_{[m]}$, and $k\leq m-1$, then
\[\Pi^{(n),m}(t)=\mathrm{Frag}(\Pi^{(n),m}(t-),\pi^{f}_{|[m]},k).\]
\end{itemize}
We verify now the compatibility property of the processes $(\Pi^{(n),m}(\cdot),m\geq 1)$ for fixed $n$.
\begin{lemma}\label{compatibility} Let $n\in \mathbb{N}$. For any $m\geq 1$ and any $t\geq 0$, $\Pi^{(n),m+1}(t)_{|[m]}=\Pi^{(n),m}(t)$ a.s.
\end{lemma}
\begin{proof} Let $m\geq 1$ and $(t_i^{c,(m+1)},i\geq 1)$ be the atoms of time of $\mathrm{PPP}_C$ whose partitions verify $\pi^{c}_{|[m+1]}\neq 0_{[m+1]}$. Similarly, denote by $(t_i^{f,(m+1)},i\geq 1)$ the atoms of time of $\mathrm{PPP}_F$ whose partitions are such that $\pi^{f}_{|[m+1]}\neq 1_{[m+1]}$ and whose index $k$ satisfies $k\leq m$. Recall that by convention for any partition $\pi$, if $\#\pi<\infty$ and $i>\#\pi$ then $\pi_i=\emptyset$.

Assume first that $t_1^{f,(m+1)}<t_1^{c,(m+1)}$. When $t<t_1^{f,(m+1)}$, we have  $\Pi^{(n),m+1}_{|[m]}(t-)=\Pi^{(n),m+1}_{|[m]}(0)=\Pi^{(n),m}(0)$ and if $t=t_1^{f,(m+1)}$ then
\begin{align*}
\Pi^{(n),m+1}(t)_{|[m]}&=\mathrm{Frag}(\Pi^{(n),m+1}(t-),\pi^f,k)_{|[m]}\\
&=\{\Pi_i^{(n),m+1}(t-)\cap [m], i\neq k, \Pi_k^{(n),m+1}(t-)\cap \pi_j^{f}\cap [m], j\geq 1\}^{\downarrow}\\
&=\{\Pi_i^{(n),m}(t-), i\neq k, \Pi_k^{(n),m}(t-)\cap \pi_j^{f}, j\geq 1\}^{\downarrow}\\
&=\mathrm{Frag}(\Pi^{(n),m}(t-),\pi^f,k)=\Pi^{(n),m}(t).
\end{align*}
Let $t_1^{f,(m+1)}<t<t_1^{c,(m+1)}$, then $\Pi^{(n),m+1}(t)_{|[m]}=\Pi^{(n),m+1}(t_1^{f,(m+1)})_{|[m]}=\Pi^{(n),m}(t)$. Denote by $j$ the index such that $j\leq m$ and $\#(\pi^{c}_j\cap [m+1])\geq 2$. 
We have for $t=t_1^{c,(m+1)}$
\begin{align*}
\Pi^{(n),m+1}(t)_{|[m]}&=\mathrm{Coag}(\Pi^{(n),m+1}(t-),\pi^c)_{|[m]}\\
&=\{\Pi_i^{(n),m+1}(t-)\cap [m], i\neq j, \cup_{i\in \pi^{c}_j\cap[m+1]}\Pi_i^{(n),m+1}(t-)\cap [m]\}\\
&=\{\Pi_i^{(n),m}(t-), i\neq k, \cup_{i\in \pi^{c}_j\cap[m]}\Pi_i^{(n),m}(t-)\}\\
&=\mathrm{Coag}(\Pi^{(n),m}(t-),\pi^{c}_{|[m]})=\Pi^{(n),m}(t)
\end{align*}
where the third equality holds since by the ordering of the subsets $\Pi_{m+1}^{(n),m+1}\cap [m]=\emptyset$.
The case $t_1^{c,(m+1)}<t_1^{f,(m+1)}$ is treated similarly and by induction $\Pi^{(n),m+1}(t)_{|[m]}=\Pi^{(n),m}(t)$ holds for any time $t\geq 0$ a.s. 
\end{proof}

The compatibility property established in Lemma \ref{compatibility} allows us to construct an $\mathcal{S}$-valued process $(\Pi^{(n)}(t),t\geq 0)$  by setting 
$\Pi^{(n)}(t):=\cup_{m\geq 1}\Pi^{(n),m}(t)$ for all $t\geq 0$. Note that by definition, when $n=\infty$, $(\Pi^{(\infty)}(t),t\geq 0)=(\Pi(t),t\geq 0)$ a.s. 

\begin{lemma}\label{coupling1} Assume that the initial partition $\Pi(0)$ has blocks with infinite sizes. For any $n\in \mathbb{N}$, set $(N^{(n)}_t,t\geq 0):=(\#\Pi^{(n)}(t),t\geq 0)$ and $\zeta^{(n)}:=\inf\{t>0; N_{t-}^{(n)}\ \mathrm{or}\ N_t^{(n)}=\infty\}$. The process $(N^{(n)}_t,0\leq t< \zeta^{(n)})$ has the same law as $(\#\Pi(t),0\leq t< \zeta)$ started from $n$. Moreover almost surely, for all $n\in \mathbb{N}$ and all $t\geq 0$, $N^{(n+1)}_t\geq N^{(n)}_t$ and $\underset{n\rightarrow \infty }{\lim} N^{(n)}_t =\#\Pi(t)$. 
\end{lemma}
\begin{proof}
By the Poisson construction of $(\Pi^{(n)}(t),t\geq 0)$, at an atom $(t,\pi^{f},j)$ of $\text{PPP}_F$, if $j\leq \#\Pi^{(n)}(t-)=N_{t-}^{(n)}$, the process $(N_t^{(n)},t\geq 0)$ jumps from $N_{t-}^{(n)}$ to
\begin{equation}\label{stateafterfrag}
\#\left(\mathrm{Frag}(\Pi^{(n)}(t-),\pi^{f},j)\right)=N_{t-}^{(n)}-1+\#\{\Pi^{(n)}_j(t-)\cap \pi^{f}_{\ell}, \ell\leq \#\pi^{f}\}.
\end{equation}
Since $\Pi_1(0),\cdots, \Pi_n(0)$ are assumed to be infinite, and the fragmentation measure is supported by partitions with no singletons, the set $\Pi^{(n)}_j(t-)$ is infinite. By using the same argument as in the proof of Proposition \ref{generatorlemma}, for all $\ell\leq \#\pi^f$, $\Pi^{(n)}_j(t-)\cap \pi^{f}_\ell\neq \emptyset$ almost surely. Hence, the state after time $t-$, \eqref{stateafterfrag} is $N_{t-}^{(n)}+k$ with $k=\#\pi^{f}-1$.

Consider now $(t,\pi^{c})$ an atom of $\text{PPP}_C$, and apply the operator $\mathrm{Coag}$. Let $j$ be the index of the non-singleton block in $\pi^c$: namely $\#\pi^{c}_j\geq 2$, then, conditionally on $\{N^{(n)}_{t-}=m\}$
\begin{equation*}\label{stateaftercoag}
\#\mathrm{Coag}(\{\Pi_1^{(n)}(t-),\cdots, \Pi_m^{(n)}(t-)\},\pi^{c})=\#\{\cup_{i\in \pi_j^{c}\cap [m]}\Pi^{(n)}_{i}(t-), \Pi^{(n)}_{\ell}(t-), \ \ell \notin \pi_j^{c}\}
\end{equation*}
and we see that $N_t^{(n)}=N_{t-}^{(n)}-k+1$ with $k:=\#(\pi_j^{c}\cap [m])$. This occurs at rate $\binom{m}{k}\lambda_{m,k}$. We deduce that $(N_t^{(n)},0\leq t<\zeta^{(n)})$ is Markov and has the same dynamics as $(\#\Pi(t), 0\leq t< \zeta)$ started from $n$, stopped at its first explosion time. 

Recall $N_t^{(n)}=\#\Pi^{(n)}(t)$ for all $t\geq 0$ and $n\in \mathbb{N}$. We now show that for all $n\geq 1$ and all $t\geq 0$, $N^{(n+1)}_t\geq N^{(n)}_t$ a.s. Let $m\in \mathbb{N}$. We check first that $\#\Pi^{(n)}_{|[m]}(t)\leq \#\Pi^{(n+1)}(t)$ for any $t\geq 0$ a.s.
Let $t>0$. If $\#\Pi^{(n)}_{|[m]}(t-)\leq \#\Pi^{(n+1)}(t-)$ then 
\begin{itemize}
\item[(i)] if $(t,\pi^c)$ is an atom of $\mathrm{PPP}_C$, by applying Lemma \ref{cardcoal}, we have \[\#\Pi^{(n)}_{|[m]}(t)=\#\mathrm{Coal}(\Pi^{(n)}_{|[m]}(t-),\pi^c)\leq\#\mathrm{Coal}(\Pi^{(n+1)}(t-),\pi^c)=\#\Pi^{(n+1)}(t),\]
\item[(ii)] if $(t,\pi^f,j)$ is an atom of $\mathrm{PPP}_F$ 
and further $j\leq \#\Pi^{(n)}_{|[m]}(t-)$ 
then
\begin{align*}
\#\Pi^{(n)}_{|[m]}(t)&=\#\left(\mathrm{Frag}(\Pi^{(n)}_{|[m]}(t-),\pi^{f},j)\right)\\
&=\#\Pi^{(n)}_{|[m]}(t-)-1+\#\{\Pi_j^{(n)}(t-)\cap \pi^f_i\cap [m],i\geq 1\}\\
&\leq \#\Pi^{(n+1)}_{|[m]}(t-)-1+\#\pi^f\leq \#\Pi^{(n+1)}(t).
\end{align*}
If $j\geq \#\Pi^{(n)}_{|[m]}(t-)+1$, then
$\#\Pi^{(n)}_{|[m]}(t)=\#\Pi^{(n)}_{|[m]}(t-)\leq \#\Pi^{(n+1)}(t-)\leq \#\Pi^{(n+1)}(t).$
\end{itemize}
Clearly $\#\Pi^{(n)}_{|[m]}(0)\leq \#\Pi^{(n+1)}(0)<\infty$ a.s. By applying (i) and (ii) until  the first explosion time of $(N^{(n+1)}_t,t\geq 0)$, we get that $\#\Pi^{(n)}_{|[m]}(t)\leq \#\Pi^{(n+1)}(t)=N^{(n+1)}_t$ for all $t\leq \zeta^{(n+1)}$ a.s. Assume by contradiction that $\#\Pi^{(n+1)}(t)<\#\Pi^{(n)}_{|[m]}(t)$ for a certain $t>\zeta^{(n+1)}$. Denote by $\zeta$ the last instant $s$ prior to $t$ at which $N^{(n+1)}_{s-}=\infty$. The process $\#\Pi^{(n+1)}$ has piecewise constant paths when lying in $\mathbb{N}$, and by applying (i) and (ii), we see that necessarily for any $\zeta<s\leq t$, $\#\Pi^{(n+1)}(s)<\#\Pi^{(n)}_{|[m]}(s)\leq m$. Hence $\#\Pi^{(n+1)}(\zeta)\leq \#\Pi^{(n)}_{|[m]}(\zeta)\leq m$. Since by definition $\#\Pi^{(n+1)}(\zeta-)=\infty$, the time $\zeta$ should be a coalescence time at which infinitely many blocks coalesce into less than $m$ blocks. This leads to a contradiction since by assumption $\Lambda(\{1\})=0$ and those coalescences are not possible. Finally $\#\Pi^{(n)}_{|[m]}(t)\leq \#\Pi^{(n+1)}(t)$ for all $t\geq 0$ a.s. and since $m$ is arbitrary, we have that for all $t\geq 0$, $N^{(n)}_t:=\#\Pi^{(n)}(t)=\underset{m\rightarrow \infty}{\lim} \#\Pi^{(n)}_{|[m]}(t)\leq N_t^{(n+1)}$ almost surely. We show similarly that $N^{(n)}_t\leq \#\Pi(t)$ for all $t$ almost surely by replacing $\Pi^{(n+1)}$ by $\Pi$ in the arguments above.\\

It remains to show that $\underset{n\rightarrow \infty }{\lim} N^{(n)}_t =\#\Pi(t) \text{ a.s.}$
Let $m\in \mathbb{N}$. Choose $n$ large enough such that $[m]\subset \cup_{i=1}^{n}\Pi_i(0)$. Then, $\Pi^{(n)}_{|[m]}(0)=\Pi_{|[m]}(0)$, and we see from the Poisson construction of $(\Pi_{|[m]}(t),t\geq 0)$ that $\Pi^{(n)}_{|[m]}(t)=\Pi_{|[m]}(t)$ for all $t\geq 0$ a.s. Hence, for any $m$, $\underset{n \rightarrow \infty}{\lim} N_t^{(n)}\geq \#\Pi^{(n)}_{|[m]}(t)=\#\Pi_{|[m]}(t)$. Letting $m$ to infinity provides $\underset{n \rightarrow \infty}{\lim} N_t^{(n)}\geq \#\Pi(t)$ which allows us to conclude, since for all $t\geq 0$, $N_t^{(n)}\leq \#\Pi(t)$  a.s. 
\end{proof}
\subsection{Coming down from infinity}\label{CDIsubsec} Recall $\theta^{\star}$ defined in \eqref{theta}. In all this section, we assume that $\theta^{\star}<1$.\\

We outline here the scheme of the proof.  Denote by $\tau^{(n)}_{n_0}$ and $\zeta^{(n)}$, respectively the first passage time below $n_0$ and the first explosion time of $(N^{(n)}_t,t\geq 0)$. We obtain in Lemma \ref{boundNt}, an upper bound of the mean of $\tau^{(n)}_{n_0}\wedge \zeta^{(n)}$, which is \textit{uniform} in the initial value $n$. We shall also see in the proof of Lemma \ref{boundNt} from where the parameter $\theta^{\star}$ comes from. Next, we define in Lemma \ref{Nnm}, a sequence of processes $(\#\Pi^{m}(t),t\geq 0)_{m\geq 1}$ approaching from below $(\#\Pi(t),t\geq 0)$. Those processes are not explosive and have the same dynamics as our initial process for a certain splitting measure $\mu_m$. We establish in Lemma \ref{boundentrance}, using
the calculations in Lemma \ref{boundNt}, that these processes are coming down from infinity, and get a bound for the mean of their first passage time below a certain state $n_0$. The latter being uniform in $m$, we will be able to conclude that the process $(\#\Pi(t),t\geq 0)$ itself goes below the level $n_0$ a.s.
\begin{lemma}\label{boundNt} Let $n\in \mathbb{N}$ and let $\zeta^{(n)}$ be the first explosion time of $(N_t^{(n)},t\geq 0)$. For any $m\leq n$, set $\tau_{m}^{(n)}:=\inf\{t\geq 0, N_t^{(n)}\leq m\}$. Then, there exists $n_0$ such that if $n\geq n_0$ then
\begin{equation}\label{upperbound}\mathbb{E}[\tau^{(n)}_{n_{0}}\wedge \zeta^{(n)}]\leq \frac{2}{1-\theta^{\star}}\sum_{k=2}^{\infty}\frac{1}{\Phi(k)}.\end{equation}
\end{lemma}
\begin{remark} The right-hand side in \eqref{upperbound} is bounded uniformly in the initial state $n$.
\end{remark}
\begin{proof}
Recall $\mathcal{D}$ in \eqref{domain} and that we work under the assumption \eqref{CDI}. Define the function $g$ on $\bar{\mathbb{N}}$ by $g(1)=0$ and $g(n)=\sum_{j=2}^{n}\frac{1}{\Phi(j)}$ when $2\leq n\leq \infty$. Note that $g$ is bounded and thus belongs to $\mathcal{D}$. Moreover $g(n)\underset{n\rightarrow \infty}{\longrightarrow} g(\infty):=\sum_{j=2}^{\infty}\frac{1}{\Phi(j)}<\infty$. On the one hand we have for any $n\geq 2$, $g(n-k+1)-g(n)=-\sum_{j=n-k+2}^{n}\frac{1}{\Phi(j)}$, and since $\Phi$ is non-decreasing, for all $2\leq j\leq n$, $1/\Phi(j)\geq 1/\Phi(n)$. Therefore
\begin{equation}\label{coalg} g(n-k+1)-g(n)\leq -\frac{k-1}{\Phi(n)}.
\end{equation}
On the other hand we have for all $n\geq 1$ and $k\in \mathbb{N}$,
\begin{equation}\label{branchingg}
g(n+k)-g(n)=\sum_{j=n+1}^{n+k}\frac{1}{\Phi(j)}=\sum_{j=n+1}^{\infty} \mathbbm{1}_{\{j\leq n+k\}}\frac{1}{\Phi(j)} \text{ and } g(\infty)-g(n)=\sum_{j=n+1}^{\infty}\frac{1}{\Phi(j)}.
\end{equation}
Plugging \eqref{coalg} and  \eqref{branchingg}  in the generator $\mathcal{L}:=\mathcal{L}^c+\mathcal{L}^f$ defined in Proposition \ref{generatorlemma} yields  \begin{align}\label{Lg}
\mathcal{L}g(n)&\leq -\frac{1}{\Phi(n)}\underbrace{\sum_{k=2}^{n}\binom{n}{k}\lambda_{n,k}(k-1)}_{=\Phi(n)}+ n\sum_{j=n+1}^{\infty}\frac{1}{\Phi(j)}\left(\sum_{k=j-n}^{\infty}\mu(k)+\mu(\infty)\right).
\end{align}
Hence, setting for any $k\in \mathbb{N}$, $\bar{\mu}(k):=\mu(\{k,k+1,\cdots, \infty\})$, one has for all $n\in \mathbb{N}$
\begin{equation} \mathcal{L}g(n)\leq -1+\sum_{k=1}^{\infty}\frac{n\bar{\mu}(k)}{\Phi(k+n)}. 
\end{equation}
By assumption, $\theta^{\star}:=\underset{n\rightarrow \infty}{\limsup} \sum_{k=1}^{\infty}\frac{n\bar{\mu}(k)}{\Phi(k+n)}<1$. Let $\epsilon=\frac{1-\theta^{\star}}{2}>0$. There exists a large enough integer $n_{0}$ such that for all $n\geq n_{0}$, $\sum_{k=1}^{\infty}\frac{n\bar{\mu}(k)}{\Phi(k+n)}\leq \theta^{\star}+\epsilon=\frac{\theta^{\star}+1}{2}$ and therefore $\mathcal{L}g(n)\leq -1+\frac{\theta^{\star}+1}{2}=\frac{\theta^{\star}-1}{2}<0$. 

For any $N>n$, let $\tau_N^{+}:=\inf\{t\geq 0; N_t^{(n)}>N\}$, the stopped process $(N_{t\wedge \tau^+_N},t\geq 0)$ has generator  $\mathcal{L}^{N}g(n):=\mathcal{L}g(n)\mathbbm{1}_{\{n\leq N\}}$. Since $g$ and $\mathcal{L}^{N}g$ are bounded, by Dynkin's formula for continuous-time Markov chains, for any fixed $k>0$ and any $n\geq n_0$, 
\begin{align*}
\mathbb{E}[g(N^{(n)}_{\tau^{(n)}_{n_{0}}\wedge k \wedge \tau_N^{+}})]&=g(n)+\mathbb{E}\left[\int_{0}^{\tau^{(n)}_{n_{0}}\wedge k \wedge \tau_N^{+}}\mathcal{L}g(N_s^{(n)})\ddr s\right]\\
&\leq g(n)+\frac{\theta^{\star}-1}{2}\mathbb{E}[\tau^{(n)}_{n_{0}}\wedge k\wedge \tau_N^{+}].
\end{align*}
%
Hence
\begin{equation}\label{bound}\mathbb{E}[\tau^{(n)}_{n_{0}}\wedge k\wedge \tau_N^{+}] 
\leq \frac{2}{1-\theta^{\star}} \left(g(n)-\mathbb{E}[g(N^{(n)}_{\tau^{(n)}_{n_{0}}\wedge  k\wedge \tau_N^{+}})] \right)\leq \frac{2}{1-\theta^{\star}} g(n).
\end{equation} 
For any $n\geq n_0$, since $\tau_N^{+}$ increases towards the explosion time of $(N_t^{(n)},t\geq 0)$, $\zeta^{(n)}$ as $N$ goes to $\infty$ almost surely, we obtain by letting $k$ to $\infty$ and $N$ to $\infty$  in \eqref{bound}
\begin{equation}\label{uniformbound} \mathbb{E}[\tau^{(n)}_{n_{0}}\wedge \zeta^{(n)}]\leq \frac{2}{1-\theta^{\star}} \sum_{k=2}^{n}\frac{1}{\Phi(k)}\leq \frac{2}{1-\theta^{\star}} \sum_{k=2}^{\infty}\frac{1}{\Phi(k)}<\infty. 
\end{equation}
\end{proof} 


We now build a monotone coupling on the space of partitions.  The main idea is to introduce a partition-valued process $(\Pi^{m}(t),t\geq 0)$, in which every fragmentations creating more than $m+1$ new blocks in the process $(\Pi(t),t\geq 0)$, are  creating at most $m$ new blocks in $(\Pi^{m}(t),t\geq 0)$. For any $m\in \mathbb{N}$, define the map \[r_m: \pi\mapsto (\pi_1,...,\pi_{m},\cup_{i=m+1}^{\infty}\pi_i).\]
By definition $r_m$ maps $\mathcal{P}_\infty$ into partitions with at most $m+1$ blocks. Set $\mu_{\mathrm{Frag}}^{m}:=\mu_{\mathrm{Frag}}\circ r_m^{-1}$. 

Let $n\geq 1$. We call respectively $(\Pi^{m}(t),t\geq 0)$ and  $(\Pi^{m, (n)}(t),t\geq 0)$, the $\mathcal{P}_\infty$-valued process, started from $\Pi(0)$, and the $\mathcal{S}$-valued Markov process, started from $(\Pi_1(0),\cdots, \Pi_n(0))$, that are constructed in a Poisson way, as $(\Pi(t),t\geq 0)$ and $(\Pi^{(n)}(t),t\geq 0)$, but with $\text{PPP}_C$ and the image of $\text{PPP}_F$ by $r_m$.
The hypothesis  $\theta^{\star}<1$ is not needed for the next two Lemmas \ref{jumpsimul} and \ref{Nnm} to hold true. They are included in this section as they will be used only for the coming down from infinity. 
\begin{lemma}\label{jumpsimul} For any $m\geq 1$, $(\Pi^{m}(t),t\geq 0)$ and $(\Pi(t),t\geq 0)$ jump simultaneously.
\end{lemma}
\begin{proof}
By construction, the atoms of coalescence are exactly those of $\text{PPP}_C$ and those of fragmentation are the images of the atoms of $\text{PPP}_F$ by $r_m$, that is to say, 
\begin{equation}\label{rm}
r_m(\pi^{f})_{|[n]}=(\pi^{f}_1\cap [n],...,\pi^{f}_{m}\cap [n],\cup_{i=m+1}^{\infty}\pi^{f}_i\cap [n]), \text{ for any }n\in \bar{\mathbb{N}}.
\end{equation}
On the one hand, if $\#\pi^{f}\leq m$ then $r_{m}(\pi^{f})=\pi^{f}$ and $\#r_{m}(\pi^{f})=\#\pi^{f}$. On the other, if $\#\pi^{f}\geq m+1$, then $\#r_{m}(\pi^{f})=m+1$. One also easily checks from \eqref{rm} that for any $m\in \mathbb{N}$ and any $n\in \mathbb{N}$, $r_{m}(\pi^{f})_{|[n]}=1_{[n]}$ if and only if $\pi^{f}_{|[n]}=1_{[n]}$. Therefore, the processes $(\Pi^{m}(t),t\geq 0)$ and $(\Pi(t),t\geq 0)$ jump simultaneously. 
\end{proof}

Recall $\bar{\mu}(m)=\mu(\{m,\cdots, \infty\})$ and denote by $\mathcal{L}^{c}$  the coalescent part of the generator $\mathcal{L}$ defined in \eqref{coalpart}.
\begin{lemma}\label{Nnm} For any $n\in \bar{\mathbb{N}}$ and $m\geq 1$, set $(N^{(n)}_{m}(t),t\geq 0):=(\#\Pi^{m, (n)}(t),t\geq 0)$. The process $(N^{(n)}_{m}(t),t\geq 0)$ is a non-explosive Markov process started from $n$ and has for generator, the operator $\mathcal{L}^{m}$ acting on any function $g: \mathbb{N}\rightarrow \mathbb{R}_+$ as follows 
\[\mathcal{L}^{m}g(\ell):=\mathcal{L}^{c}g(\ell)+\ell\sum_{k=1}^{m}\mu_{m}(k)(g(\ell+k)-g(\ell)) \text{, for all } \ell \in \mathbb{N},\]
where $\mu_{m}(k):=\mu(k)$ if $k\leq m-1$ and $\mu_m(m):=\bar{\mu}(m)$.
Moreover, almost surely for any $n,m\in \mathbb{N}$ and all $t\geq 0$, $N^{(n)}_{m}(t)\leq N^{(n+1)}_{m}(t)$, $\#\Pi^{m}(t)\leq \#\Pi^{m+1}(t)$ and  $$\underset{m\rightarrow \infty}{\lim} \#\Pi^{m}(t)=\#\Pi(t).$$ 
\end{lemma}
\begin{remark} 
The process $(N_m^{(\infty)}(t),t\geq 0):=(\#\Pi^{m}(t),t\geq 0)$ does not explode and has the same law (since it has the same generator)  as the block-counting process of any simple EFC process whose coagulation measure is $\Lambda$ and whose dislocation measure satisfies $\nu_{\mathrm{Disl}}(\mathcal{P}_{\mathrm{m}}^{k})=\mu_m(k)$ for all $k\in [m]$.
\end{remark}
\begin{proof}
Since by assumption $\mu_{\mathrm{Frag}}$ is supported by partitions whose blocks have infinite size, the blocks of any atom $\pi^{f}$ of $\text{PPP}_F$ are of infinite size, and the partition $r_m(\pi^{f})$ has thus also blocks of infinite size. Similarly as in Lemma \ref{coupling1}, replacing $\pi^{f}$ by $r_m(\pi^{f})$ in \eqref{stateafterfrag}, this guarantees that the process $(N^{(n)}_{m}(t),t\geq 0)$ is Markov. One also plainly checks that it has the same negative jumps rates as $(N^{(n)}_t,t\geq 0)$. At any fragmentation event, the block of $\Pi^{m}$ that is involved, can be splitted at most into $m+1$ sub-blocks. Therefore the positive jumps are driven by the measure $\mu_{m}$ defined over $[|1,m|]$ by $\mu_{m}(k):=\mu(k)$ if $k\leq m-1$ and $\mu_{m}(m):=\bar{\mu}(m)$.  In particular, since the process
$(N^{(n)}_{m}(t),t\geq 0)$  stays below a discrete branching process whose reproduction measure $\mu_m$ has finite support, it cannot explode. Lemma \ref{cardcoal} entails that for any fixed $m\in \mathbb{N}$, $N_m^{(n+1)}(t)\geq N_m^{(n)}(t)$ for any $t\geq 0$ and any $n\in \mathbb{N}$ a.s. We now justify that for all $n\in \mathbb{N}$ and all $m\in \mathbb{N}$, \begin{equation}\label{monotonicity}
N^{(n)}_{m}(t)\leq N^{(n)}_{m+1}(t),  \text{ for all } t\geq 0.
\end{equation}
Both processes start from $n$, and by Lemma \ref{jumpsimul} make a positive jump at the same atoms of time of $\text{PPP}_F$. Let $t$ be such an atom of time. The jump of $N^{(n)}_{m+1}$ at time $t$, is of size at most $m+1$, whereas that of $N^{(n)}_{m}$ is of size at most $m$. On the other hand, at any atom of coalescence $(t,\pi^c)$, if $N^{(n)}_{m}(t-)\leq N^{(n)}_{m+1}(t-)$, then by Lemma \ref{cardcoal} \[\#\mathrm{Coag}(\Pi^{m,(n)}(t-),\pi^c)\leq \#\mathrm{Coag}(\Pi^{m+1,(n)}(t-),\pi^c)\]
and $N^{(n)}_{m}(t)\leq N^{(n)}_{m+1}(t)$. At all jumps, the order is preserved and \eqref{monotonicity} is true for all $t$ almost surely. One can check, similarly as in the proof of Lemma \ref{coupling1}, that $N_m^{(n)}(t)$ increases  towards $\#\Pi^{m}(t)$ for any $t\geq 0$, as $n$ goes to $\infty$ almost surely. Letting $n$ to $\infty$, in \eqref{monotonicity} provides also $\#\Pi^{m}(t)\leq \#\Pi^{m+1}(t)$ for any $t\geq 0$ and any $m\geq 1$.

Last, we show now that $\underset{m\rightarrow \infty}{\lim} \#\Pi^{m}(t)=\#\Pi(t)$. Plainly, by construction for any $t\geq 0$, $\#\Pi^{m}(t)\leq\#\Pi(t)$ a.s.  By definition of the map $r_m$, it can be checked that for any partition $\pi$, and any $n\geq 1$, if $m\geq\#\pi_{|[n]}$ then $r_m(\pi)_{|[n]}=\pi_{|[n]}$. Since, there are only finitely many atoms of $\text{PPP}_C$ and $\text{PPP}_F$ on the interval of time $[0,t]$ that are seen by the process $\Pi^{m}_{|[n]}$, one can define
 \[m_n(t):=\max \{\#\pi^{f}_{|[n]}:\  \pi^{f} \text{ atom  of } \text{PPP}_F \text{ in }[0,t] \text{ such that } \pi^{f}_{|[n]}\neq 1_{[n]}\}<\infty.\]
By construction, for any $t\geq 0$,
$\Pi^{m_{n}(t)}_{|[n]}(t)=\Pi_{|[n]}(t)$ almost surely and thus $$\#\Pi_{|[n]}(t) =\#\Pi^{m_n(t)}_{|[n]}(t)\leq \#\Pi^{m_n(t)}(t) \text{ for any } t \text{ a.s.}$$ 
By monotonicity, $\#\Pi^{m_n(t)}(t)\leq \#\Pi^{\infty}(t):=\underset{m\rightarrow \infty}{\lim} \#\Pi^{m}(t)$ for any $t\geq 0$ and we have that for any $t\geq 0$
\begin{equation*} \#\Pi_{|[n]}(t)\leq \#\Pi^{\infty}(t) \text{ a.s. }
\end{equation*}
Letting $n$ to $\infty$ in the inequality above yields $\#\Pi(t)\leq \#\Pi^{\infty}(t) \text{ a.s. }$ which entails
\begin{equation}\label{sameatthelimit} \#\Pi(t)=\#\Pi^{\infty}(t) \text{ for any } t \text{ a.s. }
\end{equation}
\end{proof}
For any $m,n_0\in \mathbb{N}$, consider the first entrance times $\tau^{(n)}_{n_0,m}:=\inf\{t>0; N_m^{(n)}(t)\leq n_0\}$ and $\tau_{n_0,m}:=\inf\{t>0; \#\Pi^m(t)\leq n_0\}$. We study their limits  as $n$ and $m$ goes to infinity respectively. 
\begin{lemma}\label{conventrance} For any $n_0\in \mathbb{N}$ and any $m\in \mathbb{N}$, 
$\underset{n\rightarrow \infty}{\lim}\tau^{(n)}_{n_0,m}=\tau_{n_0,m}$ a.s. If moreover for any $m\in \mathbb{N}$ and any $s>0$, $\#\Pi^{m}(s)<\infty$, then $\underset{m\rightarrow \infty}{\lim} \tau_{n_0,m}=\tau_{n_0}$ a.s. 
\end{lemma}
\begin{proof}
Let $n_0\in \mathbb{N}$. Recall that for all $t\geq 0$, all $m\in \mathbb{N}$, $N_m^{(n)}(t)\leq N_{m}^{(n+1)}(t)$ and $\underset{n\rightarrow \infty}{\lim} N_m^{(n)}(t)=\#\Pi^{m}(t)$ a.s. Hence $\tau^{(n)}_{n_0,m}\leq \tau^{(n+1)}_{n_0,m}\leq \tau_{n_0,m}$ a.s and by letting $n$ towards infinity, we get $\tau^{(\infty)}_{n_0,m}:=\underset{n\rightarrow \infty}{\lim} \tau^{(n)}_{n_0,m} \leq \tau_{n_0,m}$ a.s. Assume by contradiction that there exists $t>0$ such that $\tau^{(\infty)}_{n_0,m}<t< \tau_{n_0,m}$. For any $n\geq 1$, $\tau^{(n)}_{n_0,m}\leq \tau^{(\infty)}_{n_0}<t$, thus there exists a time $s_n\in(0,t)$ such that $N_m^{(n)}(s_n)\leq n_0$ a.s. The sequence $(s_n)_{n\geq 1}$ is bounded by $t$, and thus converges, up to a subsequence, to some $s\in [0,t]$. Since the process $(N^{(n)}_m(u),u>0)$ lies in $\mathbb{N}$ and has piecewise constant paths, there exists $\eta_n>0$ such that for any $u\in (s-\eta_n,s+\eta_n)\cap [0,t]$, $N_m^{(n)}(u)=N_m^{(n)}(s)$. Let $(s_{\varphi(n)},n\geq 1)$ be a subsequence  such that for all $n\geq 1$, $s_{\varphi(n)}\in (s-\eta_n,s+\eta_n)\cap [0,t]$ . Then, since $\varphi(n)\geq n$, \[N_m^{(n)}(s)=N_m^{(n)}(s_{\varphi(n)})\leq N_m^{(\varphi(n))}(s_{\varphi(n)})\leq n_0.\]
Finally, $\#\Pi^{m}(s)=\underset{n\rightarrow \infty}{\lim} N_m^{(n)}(s)\leq n_0$, thus $\tau_{n_0,m}\leq s$ which contradicts the fact that $\tau_{n_0,m}>t$. Hence, $\underset{n\rightarrow \infty}{\lim} \tau^{(n)}_{n_0,m}=\tau_{n_0,m}$ a.s. 

The convergence $\underset{m\rightarrow \infty}{\lim} \tau_{n_0,m}=\tau_{n_0}$ will follow from similar arguments. Note first that $\underset{m\rightarrow \infty}{\lim} \tau_{n_0,m}\leq \tau_{n_0}$. As previously, assume that $\underset{m\rightarrow \infty}{\lim} \tau_{n_0,m}<t< \tau_{n_0}$. One can find a convergent sequence $(s_m)_{m\geq 1}$ such that $0<\tau_{n_0,2}\leq s_m<t$ for any $m\in \mathbb{N}$ and $\#\Pi^{m}(s_m)\leq n_0$. Let $s:=\underset{m\rightarrow \infty}{\lim} s_m$, we have $s\in [\tau_{n_0,2}, t]$ and by the assumption $\#\Pi^m(s)<\infty$. Therefore, there exists $\eta_m>0$ such that for all $u\in (s-\eta_m,s+\eta_m)$, $\#\Pi^{m}(u)=\#\Pi^m(s)$. By choosing a subsequence $(s_{\varphi(m)},m\geq 1)$ such that $s_{\varphi(m)}\in (s-\eta_m,s+\eta_m)\cap [\tau_{n_0,2},t]$, we see that \[\#\Pi^{m}(s)=\#\Pi^{m}(s_{\varphi(m)})\leq \#\Pi^{\varphi(m)}(s_{\varphi(m)})\leq n_0.\] Recall Lemma \ref{Nnm} and that $\#\Pi^{m}(s)\underset{m\rightarrow \infty}{\longrightarrow} \#\Pi(s)$ a.s. We conclude as before by the contradiction $\tau_{n_0}\leq s$ and $\tau_{n_0}>t$.
\end{proof}
We are now ready to finish the proof. Recall $\theta^{\star}:=\underset{n\rightarrow \infty}{\limsup}\sum_{k=1}^{\infty}\frac{n\bar{\mu}(k)}{\Phi(n+k)}$ and the assumption $\theta^{\star}<1$. 
\begin{lemma}\label{boundentrance} There exists a large enough
integer $n_0$ such that 
$$\mathbb{E}(\tau_{n_0})\leq \frac{2}{1-\theta^{\star}}\sum_{k=2}^{\infty}\frac{1}{\Phi(k)}<\infty.$$
\end{lemma}
\begin{proof}
As in the proof of Lemma \ref{boundNt}, consider $n_0$ large enough such that for all $n\geq n_0$, $\sum_{k=1}^{\infty}\frac{n\bar{\mu}(k)}{\Phi(k+n)}\leq \theta^{\star}+\frac{1-\theta^{\star}}{2}=\frac{\theta^{\star}+1}{2}.$
According to Lemma \ref{Nnm}, $(N^{(n)}_{m}(t),t\geq 0)$ has for generator $\mathcal{L}^{m}$.  Equation \eqref{Lg} applied to the process $(N^{(n)}_{m}(t),t\geq 0)$ gives for all $n\geq n_0$
\[\mathcal{L}^{m}g(n)\leq -1+\sum_{k=1}^{\infty}\frac{n\bar{\mu}_{m}(k)}{\Phi(n+k)}=-1+\sum_{k=1}^{m}\frac{n\bar{\mu}(k)}{\Phi(n+k)}\leq -1+\sum_{k=1}^{\infty}\frac{n\bar{\mu}(k)}{\Phi(n+k)}\leq \frac{\theta^{\star}-1}{2}<0.\]
Therefore, for any $m\geq 1$ and $n\geq n_0$, $\mathbb{E}[\tau^{(n)}_{n_0,m}\wedge \zeta^{(n)}_{m}]\leq \frac{2}{1-\theta^{\star}}\sum_{k=2}^{\infty}\frac{1}{\Phi(k)}$
with $\zeta^{(n)}_{m}:=\inf \{t>0; N_m^{(n)}(t-)=\infty\}$. By Lemma \eqref{Nnm}, the process $(N_m^{(n)}(t),t\geq 0)$ does not explode and therefore  $\zeta^{(n)}_{m}=\infty$ a.s. Hence, we get
$$\mathbb{E}[\tau^{(n)}_{n_0,m}]\leq \frac{2}{1-\theta^{\star}}\sum_{k=2}^{\infty}\frac{1}{\Phi(k)}.$$
By Lemma \ref{conventrance}, $\tau^{(n)}_{n_0,m}$ increases towards  $\tau_{n_0,m}:=\inf\{t\geq 0, \#\Pi^{m}(t)\leq n_0\}$ as $n$ goes to $\infty$. By monotone convergence, we see that for any $m\geq 1$, 
\begin{equation}\label{boundtaum}
\mathbb{E}[\tau_{n_0,m}]\leq \frac{2}{1-\theta^{\star}}\sum_{k=2}^{\infty}\frac{1}{\Phi(k)}.\end{equation}
Therefore the process $(\#\Pi^{m}(t),t\geq 0)$ comes down from infinity and since it does not explode, we have that $\#\Pi^{m}(s)<\infty$ for any $m\geq 1$ and all $s>0$ a.s. 
%
By Lemma \ref{conventrance}, $\underset{m\rightarrow \infty}{\lim} \tau_{n_0,m}=\tau_{n_0}$ a.s. and we obtain  by letting $m$ to $\infty$ in \eqref{boundtaum} 
$$\mathbb{E}[\tau_{n_0}]\leq \frac{2}{1-\theta^{\star}}\sum_{k=2}^{\infty}\frac{1}{\Phi(k)},$$
where we recall $\tau_{n_0}=\inf\{t\geq 0, \#\Pi(t)\leq n_0\}$. This achieves the proof.
\end{proof}
\begin{remark}\label{remarkoncoupling}
If one drops the assumption that the fragmentation measure is supported by partitions with no singleton blocks, then the process $(\#\Pi^{m}(t),t\geq 0)$, defined in Lemma \ref{Nnm}, is not Markov. Indeed, at an atom $(t,\pi^{f},k)$ of fragmentation, the number of blocks in $\Pi^{m}$ evolves as follows
$$\#\Pi^{m}(t)-\#\Pi^{m}(t-)=-1+\#\{\Pi^{m}_k(t-)\cap r_m(\pi^{f})_j, 1\leq j\leq \#r_m(\pi^{f})\}.$$
If $\pi^{f}$ has singletons then the partition $r_m(\pi^{f})$ would have (finitely many) singletons with positive probability. Thus, on the event $\{r_m(\pi^{f})_j=\{i\} \text{ and } i\notin \Pi^{m}_k(t-)\}$, the set $\Pi^{m}_k(t-)\cap r_m(\pi^{f})_j$ is empty and the jump size of $\#\Pi^{m}$ is not $\#r_m(\pi^f)$ but depends on the constituent elements of the blocks of $\Pi^{m}(t-)$. 
\end{remark}

\begin{remark} The arguments involving the non-explosive processes $(N^{(n)}_m(t),t\geq 0, n\geq 1)$, approaching from below $(N^{(n)}_t,t\geq 0)$, in a  monotone way, are reminiscent to those used  in \cite[Section 7]{FoucartEJP} for constructing logistic continuous-state branching processes reflected at $\infty$.
\end{remark}

\subsection{Staying infinite}\label{Stayinfinitesubsec} Recall $\theta_{\star}$ in \eqref{theta} and the assumption \eqref{CDI}: $\sum_{k=2}^{\infty}\frac{1}{\Phi(k)}<\infty$. In all this section, we assume that $\theta_{\star}>1$. We shall establish the second part of Theorem \ref{mainthm}, namely that the process $(\Pi(t),t\geq 0)$ stays infinite. We argue by contradiction and assume from now on that the process does come down from infinity.

We outline here the scheme of the proof. We shall see that when the process comes down from infinity, the first jump time at which the process loses a proportion $p\in (0,1)$ of its blocks is strictly positive a.s. (Lemma \ref{controljumps}). Next, we make use of the function $f(n):=\sum_{j=n+1}^{\infty}\frac{1}{\Phi(j)}$, and find a martingale argument entailing that before this jump time, the process has actually infinitely many blocks (Lemmas \ref{beforep}  and \ref{stayinfinite}). The contradiction will lie on the fact that the coming down from infinity is instantaneous (Lemma \ref{zeroone}).\\ 

We need first the following lemmas (lifted from \cite[Lemmas 6.2 and 6.3]{coaldist11}).
\begin{lemma}\label{concentration} For any $p\in(0,1)$. There exists $x_p\in (0,1)$ such that if $x\in (0,x_p)$ and $(X_k,k\geq 1)$ is a sequence of i.i.d Bernoulli random variables with parameter $x$ then for any $n_0\geq 1$, there is a positive constant $C_{p,n_0}$  such that
\[\mathbb{P}\left(\text{there exists } n\geq n_0, \sum_{k=1}^{n}X_k\geq np\right)\leq C_{p,n_0}x^{n_0 p}.\]

\end{lemma}
\begin{proof}
By the Markov inequality, \[\mathbb{P}\left(\sum_{k=1}^{n}X_k\geq np\right)\leq e^{-npt}\mathbb{E}[e^{t\sum_{k=1}^{n}X_k}]=e^{-n(pt-\log(e^{t}x+1-x))}.\]
When choosing $t=\log(1/x)$, we get the bound $\mathbb{P}\left(\sum_{k=1}^{n}X_k\geq np\right)\leq e^{-nh(x)}$ with \[h(x):=p\log(1/x)-\log(2-x).\] In particular, since $h(x)\underset{x\rightarrow 0}{\longrightarrow} \infty$,  there exists $x_p\in (0,1)$ such that for any $x\in (0,x_p)$, $h(x)>0$ and
we get \[\mathbb{P}\left( \exists n\geq n_0, \sum_{k=1}^{n}X_k\geq np\right)\leq \frac{e^{-n_0 h(x)}}{1-e^{-h(x)}}\leq C_{p,n_0}x^{n_0 p}\]
with $C_{p,n_0}=2^{n_0}\sup_{x\in (0,x_p)}1/(1-e^{-h(x)})\in (0,\infty).$
\end{proof} 
\begin{lemma}\label{controljumps}
Assume that the process $(\Pi(t),t\geq 0)$  comes down from infinity. For any $p\in (0,1)$, the first jump which makes decrease $(\#\Pi(t),t\geq 0)$ by a proportion of size at least $p$ is strictly positive a.s. Namely \[\sigma_p:=\inf\{t>0; \ \#\Pi(t)\leq (1-p) \#\Pi(t-)\}>0 \text{ a.s.}\]
Moreover, setting
\[\sigma_{p}^{(n)}:=\inf\{t\geq 0; N_t^{(n)}\leq (1-p)N_{t-}^{(n)}\},\]
we have that $\sigma_p^{(n)}\underset{n\rightarrow\infty}{\longrightarrow} \sigma_{p}$ a.s.
\end{lemma}

\begin{proof} Obviously, only coalescence times can make decrease the number of blocks. Since the process $\Pi$ is c\`adl\`ag, $\Pi(0)=\Pi(0+)$ and $0$ is not a jump time of $(\#\Pi(t),t\geq 0)$. It remains to explain why $\sigma_p$ is not an accumulation of coalescence times near $0$.  The jumps that make decrease the number of blocks by a fraction $p$ are atoms $(t,\pi^{c})$ of $\text{PPP}_C$ satisfying $\#\Pi(t-)<\infty$ and \begin{equation}\label{losefractionp}\#\Pi(t)=\#\mathrm{Coag}(\Pi(t-),\pi^{c})\leq \#\Pi(t-)(1-p).\end{equation}
Recall the i.i.d random variables $(X_k,k\geq 1)$ defined in Section \ref{cdicoal}. By definition, for any $n\in \mathbb{N}$, $\sum_{k=1}^{n}X_k$ equals $\#(\pi_i^{c}\cap [n])$, where $\pi_i^{c}$ is the non-singleton block of $\pi^c$. We see that jumps satisfying \eqref{losefractionp} occurring before $1$ and $\tau_{n_0}$ are elements of
\[J_{p}:=\left\{(t,\pi^{c}); t\leq 1 \text{ and } \exists n\geq n_0; \sum_{k=1}^{n}X_k\geq np\right\}.\]
Recall  Lemma \ref{concentration} and choose $n_0\geq 2/p$. By the compensation formula of Poisson point process 
\begin{align*}
\mathbb{E}\big(\text{PPP}_C(J_p)\big)&=\int_{0}^{1}\mathbb{P}\big(\exists n\geq n_0; \sum_{k=1}^{n}X_k\geq np\big )\nu_{\mathrm{Coag}}(\ddr x)\\
&\leq C_{p,n_0}\int_{0}^{x_p}x^{n_0 p}\nu_{\mathrm{Coag}}(\ddr x)+\int_{x_p}^{1}\nu_{\mathrm{Coag}}(\ddr x)\\
&\leq C_{p,n_0}\int_{0}^{x_p}x^{2}\nu_{\mathrm{Coag}}(\ddr x)+\frac{1}{x_p^{2}}\int_{x_p}^{1}x^{2}\nu_{\mathrm{Coag}}(\ddr x)<\infty.
\end{align*}
Finally, $\text{PPP}_C(J_p)<\infty$ a.s. and  there is only a finite number of jumps satisfying \eqref{losefractionp} before $\tau_{n_0}$. Since $0$ is not one of them, we have $\sigma_p \wedge \tau_{n_0}>0$ a.s. which entails $\sigma_p>0$ a.s. 

As by the assumption \eqref{simpleEFC}, there are no coagulations of all blocks at once, (namely $\Lambda$ has no mass at $1$), necessarily $\#\Pi(\sigma_{p}-)<\infty$ a.s. Therefore, there exists $\epsilon>0$  such that for all $u\in (\sigma_p-\epsilon,\sigma_p)$, $\#\Pi(u)=\#\Pi(\sigma_p-)$ and for all $u\in (\sigma_p,\sigma_p+\epsilon)$, $\#\Pi(u)=\#\Pi(\sigma_p)<\infty$. Recall Lemma \ref{coupling1}. As $N_t^{(n)}$ increases towards $\#\Pi(t)$ a.s when $n$ goes to $\infty$, there is a large enough $n_0$ such that, for all $n\geq n_0$ and all $t\in (\sigma_p-\epsilon,\sigma_p+\epsilon)$,  $N_t^{(n)}=\#\Pi(t)$. Thus, $\sigma_p=\sigma^{(n)}_p$ for all $n\geq n_0$ and the last convergence statement is established.  
\end{proof}
\begin{lemma}\label{convexity} For any $p\in (0,1)$, for large enough $x$, \[\frac{\Phi(x)}{\Phi((1-p)x)}\leq \left(\frac{1}{1-p}\right)^3.\]
\end{lemma}
\begin{proof} Recall $\Psi$ defined in \eqref{psi}. The function $\varphi:x\mapsto \Psi(x)/x$ is the Laplace exponent of driftless subordinator and is therefore a concave function satisfying $\varphi(0)=0$. Therefore
\[\varphi((1-p)x+p.0)=\frac{\Psi((1-p)x)}{(1-p)x}\geq (1-p)\varphi(x)+p\varphi(0)=(1-p)\frac{\Psi(x)}{x}.\]
Thus $\frac{\Psi(x)}{\Psi((1-p)x)}\leq \left(\frac{1}{1-p}\right)^2$ and
\begin{equation}\label{unimportantbadbound}\frac{\Psi(x)}{\Psi((1-p)x)}\frac{\Phi((1-p)x)}{\Phi(x)}\leq \frac{1}{(1-p)^2}\frac{\Phi((1-p)x)}{\Phi(x)}.\end{equation}
Recall that $\Psi(x)\underset{x\rightarrow \infty}{\sim} \Phi(x)$. Therefore, the left hand side in \eqref{unimportantbadbound} goes to $1$ as $x$ goes to $\infty$, and for large enough $x$, 
\begin{equation*}1-p\leq \frac{1}{(1-p)^2}\frac{\Phi((1-p)x)}{\Phi(x)}.
\end{equation*}
This enables us to conclude. 
\end{proof}
For any $n\in \bar{\mathbb{N}}$, define the process $(N^{(n),p}_t,t\geq 0):=(N^{(n)}_{t\wedge \sigma^{(n)}_p},t\geq 0)$. Recall the assumption \eqref{CDI} and set $f(n):=\sum_{j=n+1}^{\infty}\frac{1}{\Phi(j)}$ for any $n\geq 1$, and $f(\infty)=0$. Note that $f(n)$ decreases towards $f(\infty)=0$ as $n$ goes to $\infty$.
\begin{lemma}\label{beforep}
There exists $p\in (0,1)$ and $n_0\in \mathbb{N}$ such that for all $n\geq m\geq n_0$,
\[\mathbb{E}\big(f(N^{(n),p}_{t\wedge
 \zeta^{(n)}\wedge \tau_m^{(n)}})\big)\leq f(n)\]
 where $\zeta^{(n)}$ is the first explosion time of $(N_t^{(n)},t\geq 0)$.
\end{lemma}
\begin{proof}
Let $\zeta^{(n),p}$ be the first explosion time of $(N_t^{(n),p},t\geq 0)$, the stopped process $(N^{(n),p}_t, 0\leq t\leq \zeta^{(n),p})$ is Markov and has for generator $\mathcal{L}^{p}f:=\mathcal{L}^{c,p}f+\mathcal{L}^{f}f$ with
\begin{equation*}\label{truncgen}
\mathcal{L}^{c,p}f(n)=\sum_{k=2}^{\left \lfloor{pn} \right\rfloor}\binom{n}{k}\lambda_{n,k}\left(f(n-k+1)-f(n)\right).
\end{equation*}
Notice that $f$ is bounded and thus belongs to the domain of the generator $\mathcal{L}^{p}$ (which matches with $\mathcal{D}$ in \eqref{domain}). For any $2\leq k\leq \left \lfloor{pn} \right\rfloor$ and $j\geq n-k+2$, since $\Phi$ is non-decreasing \[\Phi(j)\geq \Phi(n-k+2)\geq \Phi((1-p)n).\]
We obtain, for large enough $n$,
\begin{equation}\label{truncgen2}
\mathcal{L}^{c,p}f(n)=\sum_{k=2}^{\left \lfloor{pn} \right\rfloor}\binom{n}{k}\lambda_{n,k}\left(\sum_{j=n-k+2}^{n}\frac{1}{\Phi(j)}\right)\leq \sum_{k=2}^{\left \lfloor{pn} \right\rfloor}\binom{n}{k}\lambda_{n,k}\frac{k-1}{\Phi((1-p)n)}\leq \frac{\Phi(n)}{\Phi((1-p)n)}.
\end{equation}
Applying Lemma \ref{convexity} in the last inequality of \eqref{truncgen2}, provides that for large enough $n$,
\begin{equation}\label{boundcp}
\mathcal{L}^{c,p}f(n)\leq \frac{1}{(1-p)^3}.
\end{equation}
We now apply the second part of the generator, $\mathcal{L}^{f}$, to the map $f$. Recall $\bar{\mu}(j)=\mu(\{j,j+1,\cdots, \infty\})$ for all $j\in \mathbb{N}$. For any $n\in \mathbb{N}$, one has
\begin{align}\label{Lff}
\mathcal{L}^{f}f(n)&=\sum_{k=2}^{\infty}n\mu(k)(f(n+k)-f(n))+n\mu(\infty)(f(\infty)-f(n)) \nonumber\\
&=-\sum_{k=2}^{\infty}n\mu(k)\sum_{j=n+1}^{n+k}\frac{1}{\Phi(j)}-n\mu(\infty)\sum_{j=n+1}^{\infty}\frac{1}{\Phi(j)}\nonumber\\
&=-\sum_{k\in \bar{\mathbb{N}}}\sum_{j}\frac{n\mu(k)}{\Phi(j)}\mathbbm{1}_{\{n+1\leq j\leq n+k\}}=-\sum_{j\geq n+1}\left(\sum_{k\geq j-n \atop k\in \bar{\mathbb{N}}}n\mu(k)\right)\frac{1}{\Phi(j)}\nonumber\\
&=-\sum_{j=n+1}^{\infty}\frac{n\bar{\mu}(j-n)}{\Phi(j)}=-\sum_{j=1}^{\infty}\frac{n\bar{\mu}(j)}{\Phi(j+n)}.
\end{align}
From the last equality and the definition of $\theta_\star$ in \eqref{theta},  we see that  \[\underset{n\rightarrow \infty}{\limsup} \ \mathcal{L}^{f}f(n)= -\theta_{\star}.\] 
Recall that by assumption $\theta_{\star}>1$. Assume first $\theta_\star<\infty$. Let $\epsilon>0$ small enough such that $\theta_\star-\epsilon>1$, there is $n_0$ such that for all $n\geq n_0$, \[\mathcal{L}^{p}f(n)=\mathcal{L}^{c,p}f(n)+\mathcal{L}^{f}f(n)\leq \frac{1}{(1-p)^3}-\theta_{\star}+\epsilon.\] 
Since $\frac{1}{(1-p)^3}\underset{p\rightarrow 0+}{\longrightarrow} 1$, one can choose a small enough $p\in (0,1)$ such that $\frac{1}{(1-p)^3}\leq \theta_{\star}-\epsilon$.  Finally one gets for all $n\geq n_0$
\begin{equation}\label{Lp} \mathcal{L}^{p}f(n)\leq 0.
\end{equation}
Plainly, when $\theta_\star=\infty$, the inequality \eqref{Lp} holds also true for large enough $n$. By Dynkin's formula, for any $n\geq m\geq n_0$
\[\mathbb{E}\big(f(N^{(n),p}_{t\wedge \tau_{m}^{(n)}\wedge \zeta^{(n),p}})\big)-f(n)=\mathbb{E}\left[\int_{0}^{t\wedge \tau_{m}^{(n)}\wedge \zeta^{(n),p}}\mathcal{L}^{p}f(N^{(n),p}_{s})\ddr s\right]\leq 0\]
It remains to see that $\mathbb{E}\big(f(N^{(n),p}_{t\wedge
 \zeta^{(n),p}\wedge \tau_m^{(n)}})\big)=\mathbb{E}\big(f(N^{(n),p}_{t\wedge \zeta^{(n)}\wedge \tau_m^{(n)}})\big)$. It suffices to check that $\zeta^{(n),p}\wedge \sigma^{(n)}_p=\zeta^{(n)}\wedge \sigma^{(n)}_p$ a.s. On the one hand, on the event $\{\zeta^{(n),p}<\sigma_p^{(n)}\}$,  $\zeta^{(n)}=\zeta^{(n),p}$ a.s. thus $\zeta^{(n),p}\wedge \sigma^{(n)}_p=\zeta^{(n)}\wedge \sigma^{(n)}_p$. On the other hand, on $\{\zeta^{(n),p}>\sigma_p^{(n)}\}$, $\zeta^{(n),p}=\infty$ a.s. thus $\zeta^{(n),p}\wedge \sigma^{(n)}_p=\zeta^{(n)}\wedge \sigma^{(n)}_p$. This ends the proof.
\end{proof}

We are now able to finish the proof by finding a contradiction.
\begin{lemma}\label{stayinfinite} If $\theta_{\star}>1$, the process stays infinite.
\end{lemma}
\begin{proof}
Recall that we assume that $\Pi(0)$ has infinitely many blocks of infinite size.  Since the process is assumed to come down from infinity, according to Lemma \ref{zeroone}, the process $(N_t,t\geq 0):=(\#\Pi(t),t\geq 0)$ leaves infinity instantaneously. Moreover, by Proposition \ref{generatorlemma}, the process $(N_t,t\geq 0)$ is Markov when lying in $\mathbb{N}$. Consider an excursion from $\infty$ with length $\zeta$ (possibly infinite), such that $\zeta>\tau_m>t$ for some $t>0$ and $m\geq n_0$. By the Markov property at time $t$, conditionally on $N_t$, the process $(N_{t+s}, 0\leq s\leq \zeta-t)$ has the same law as the process started from $N_t$ and stopped at its first explosion time. According to Lemma \ref{coupling1}, the latter has the same law as $(N^{(N_t)}_s, s\leq \zeta^{(N_t)})$ and by applying Lemma \ref{beforep}, we get
\[\mathbb{E}[f(N_{(t+s)\wedge \tau_m\wedge \zeta \wedge \sigma_p})\mathbbm{1}_{s+t<\tau_m<\zeta}]=\mathbb{E}\left(f(N^{(N_t)}_{s\wedge \tau_m\wedge \sigma_p}) \mathbbm{1}_{s<\tau_m^{(N_t)}<\zeta^{(N_t)}}\right)\leq \mathbb{E}(f(N_t)).\]
By the right-continuity of the process $(N_t,t\geq 0)$, see Proposition \ref{generatorlemma}, $N_t\underset{t\rightarrow 0+}{\longrightarrow} \infty$ a.s. Since $f$ is bounded and has limit $0$ at $\infty$, by using Lebesgue's theorem, we get that $\mathbb{E}(f(N_t))\underset{t\rightarrow 0+}{\longrightarrow}0$. Hence
$$\underset{t\rightarrow 0+}{\lim} \mathbb{E}[f(N_{(s+t)\wedge \tau_m\wedge \zeta \wedge \sigma_p})\mathbbm{1}_{s+t<\tau_m<\zeta}]=0.$$
A second application of Lebesgue's theorem yields
$$\mathbb{E}[f(N_{s\wedge \tau_m \wedge \sigma_p})\mathbbm{1}_{s\leq \tau_m<\zeta}]=0.$$
Since $f$ is positive then $f(N_{s\wedge \tau_m\wedge \sigma_p})=0$ a.s. on the event $\{s\leq \tau_m<\zeta\}$. 
This entails that if $s\leq \tau_m<\zeta$ then $N_{s\wedge \tau_m \wedge \sigma_p}=\infty$ a.s. Recall that $\sigma_p>0$ a.s. One has therefore, for $s\in (0,\sigma_p\wedge \tau_m)$, $N_{s}=\infty$ a.s, this is a contradiction since according to the zero-one law stated in Lemma \ref{zeroone}, if the process $\Pi$ does not stay infinite then it leaves $\infty$ instantaneously a.s. 
\end{proof}

We end this section by dealing with the critical boundary case $\theta_{\star}=1$ in the particular case where only binary coagulations are allowed. Kyprianou et al's result \cite[Theorem 1.1]{kyprianou2017} is thus recovered in our framework and generalized to cases where the measure $\mu$ gives mass to $\mathbb{N}$.
\begin{proposition}\label{propcritical} Let $c_{\mathrm{k}}>0$ and $\lambda>0$. Assume $\Lambda=c_{\mathrm{k}}\delta_0$ and $\mu(\infty)=\lambda$. If $\theta:=\frac{2\lambda}{c_{\mathrm{k}}}\geq 1$, then the process stays infinite. In particular, the process stays infinite in the critical case $\theta=1$.
\end{proposition}
\begin{proof}
Since $\Lambda=c_{\mathrm{k}}\delta_0$, $\Phi(k)=c_{\mathrm{k}}\binom{k}{2}$ for all $k\geq 2$.
Set  $f(n):=\sum_{k=n+1}^{\infty}\frac{1}{\Phi(k)}$ for any $n\geq 1$, one has $f(n)=\frac{2}{c_{\mathrm{k}}}\frac{1}{n}$ and $\mathcal{L}^{c}f(n)=1$ for all $n\geq 1$. Moreover, recall \eqref{Lff}, for any  $n\geq 1$, $\mathcal{L}^{f}f(n)=-\sum_{j=1}^{\infty}\frac{n\bar{\mu}(j)}{\Phi(j+n)}$ with $\bar{\mu}(j)=\mu(\{1,2,\cdots, \infty\})\geq \mu(\infty)=\lambda$. Hence for any $n\geq 1$, $\mathcal{L}^{f}f(n)\leq -\frac{2\lambda}{c_{\mathrm{k}}}n\sum_{j=n+1}^{\infty}\frac{1}{j(j-1)}=-\theta.$ Therefore, $\mathcal{L}f(n)\leq 1-\theta\leq 0$ for any $n\geq 1$ and assuming that the process comes down from infinity, the same reasoning as in the proof of Lemma \ref{stayinfinite} yields a contradiction. We conclude that the process stays infinite. 
\end{proof}
\begin{remark}\label{FEFC3}  We stress that in the proof of Proposition \ref{propcritical} the coupling between $(\Pi(t),t\geq 0)$ and $(\Pi(t\wedge \sigma_p),t\geq 0)$ is not used. 
\end{remark}
\section{Examples}\label{proofsandexamples}
We will establish in this section Corollary \ref{FEFC2}, Corollary \ref{suffcond1} and Proposition \ref{regularcdi}. We start by Corollary \ref{FEFC2} which is easily derived from Theorem \ref{mainthm}. 

A difficulty while dealing with the parameters $\theta^{\star}$ and $\theta_{\star}$, lies in the fact that the variables $n$ and $k$ are not separated in formulas \eqref{theta}. 
We give some technical lemmas providing a general recipee for studying $\theta^{\star}$ and $\theta_\star$ and decide whether it is $0$, $\infty$ or in $(0,\infty)$.
\subsection{Analysis of the parameters}

\begin{lemma}\label{boundsigma} For all $n\in \mathbb{N}$, set $\ell(n):=\sum_{k=1}^{n}\bar{\mu}(k)$. 
\begin{enumerate}
\item If $\frac{n\ell(n)}{\Phi(n)}\underset{n\rightarrow \infty}\longrightarrow \infty$, then $\theta_\star=\theta^{\star}=\infty$. 
\item Set $\underline{\theta}^{\star}:=\underset{n\rightarrow \infty}{\limsup }\frac{n \ell(n)}{\Phi(2n)} \text{ and }\overline{\theta}^{\star}:=\underset{n\rightarrow \infty}{\limsup }\frac{n \ell(n)}{\Phi(n)}$. One has 
$\theta^{\star}\geq \underline{\theta}^{\star}$ and $\frac{1}{4}\overline{\theta}^{\star}\leq \underline{\theta}^{\star}\leq \frac{1}{2}\overline{\theta}^{\star}$.
If moreover  $\underset{n\rightarrow \infty}{\limsup } \ n \sum_{k=n}^{\infty}\frac{\bar{\mu}(k)}{\Phi(k)}=0$, then $\theta^{\star}\leq \overline{\theta}^{\star}$.

\item Set $\underline{\theta}_{\star}:=\underset{n\rightarrow \infty}{\liminf }\frac{n \ell(n)}{\Phi(2n)}$ and $\overline{\theta}_{\star}:=\underset{n\rightarrow \infty}{\liminf }\frac{n \ell(n)}{\Phi(n)}$. One has $\theta_{\star}\geq \underline{\theta}_{\star}$ and $\frac{1}{4}\overline{\theta}_{\star}\leq \underline{\theta}_{\star}\leq \frac{1}{2}\overline{\theta}_{\star}$. If moreover $\underset{n\rightarrow \infty}{\liminf } \ n \sum_{k=n}^{\infty}\frac{\bar{\mu}(k)}{\Phi(k)}=0$, then $\theta_{\star}\leq  \overline{\theta}_{\star}$.
\item If $\frac{n\ell(n)}{\Phi(n)}\underset{n\rightarrow \infty}\longrightarrow 0$, then 
$\theta^{\star}=\underset{n\rightarrow \infty}{\limsup } \ n \sum_{k=n}^{\infty}\frac{\bar{\mu}(k)}{\Phi(k)}$ and $\theta_{\star}=\underset{n\rightarrow \infty}{\liminf} \ n \sum_{k=n}^{\infty}\frac{\bar{\mu}(k)}{\Phi(k)}$.
\end{enumerate}
\end{lemma}
\begin{proof}
We focus on $\theta^{\star}$. Arguments for $\theta_{\star}$ are the same replacing $\limsup$ by $\liminf$. We show that in general
\begin{equation}\label{crudebound} \underset{n\rightarrow \infty}{\limsup }\left(\frac{n \ell(n)}{\Phi(2n)}+ n \sum_{k=2n+1}^{\infty}\frac{\bar{\mu}(k)}{\Phi(k)}\right)\leq \theta^{\star}\leq \underset{n\rightarrow \infty}{\limsup } \left(\frac{n \ell(n)}{\Phi(n)}+ n \sum_{k=n}^{\infty}\frac{\bar{\mu}(k)}{\Phi(k)}\right).
\end{equation}
We have $$\sum_{k=1}^{\infty}\frac{n\bar{\mu}(k)}{\Phi(n+k)}=\sum_{k=1}^{n}\frac{n\bar{\mu}(k)}{\Phi(n+k)}+ n\sum_{k=n+1}^{\infty}\frac{\bar{\mu}(k)}{\Phi(n+k)}.$$
Since $\Phi$ is non-decreasing, one has $\frac{1}{\Phi(n+k)}\leq \frac{1}{\Phi(n)}$ for all $1\leq k\leq n$ and $\frac{1}{\Phi(n+k)}\leq \frac{1}{\Phi(k)}$ for all $k\geq n+1$.
Therefore $$\theta^{\star}\leq \underset{n\rightarrow \infty}{\limsup }\left( \frac{n \ell(n)}{\Phi(n)}+n \sum_{k=n}^{\infty}\frac{\bar{\mu}(k)}{\Phi(k)}\right).$$
On the other hand, 
\begin{align*}
\sum_{k=1}^{n}\frac{n\bar{\mu}(k)}{\Phi(n+k)}+ n\sum_{k=n+1}^{\infty}\frac{\bar{\mu}(k)}{\Phi(n+k)}
&=\sum_{k=n+1}^{2n}\frac{n\bar{\mu}(k-n)}{\Phi(k)}+n\sum_{k=2n+1}^{\infty}\frac{\bar{\mu}(k-n)}{\Phi(k)}\\
&\geq \sum_{k=1}^{n}\frac{n\bar{\mu}(k)}{\Phi(2n)}+n\sum_{k=2n+1}^{\infty}\frac{\bar{\mu}(k)}{\Phi(k)}\geq \frac{n \ell(n)}{\Phi(2n)}+n\sum_{k=2n+1}^{\infty}\frac{\bar{\mu}(k)}{\Phi(k)},
\end{align*}
and we obtain \eqref{crudebound}. As a first consequence, by replacing $\limsup$ by $\liminf$ in \eqref{crudebound}, we see that,  if $\frac{n\ell(n)}{\Phi(2n)}\underset{n\rightarrow \infty}\longrightarrow \infty$, then $\theta_\star=\infty$. Recall $\Psi$ in \eqref{psi} and that $\Psi(n)\underset{n\rightarrow \infty}{\sim}\Phi(n)$. Since $\varphi(x):=\Psi(x)/x$ is concave, for all $x\geq 0$,  $\varphi(x/2)\geq \varphi(x)/2$ and  we obtain with $x=2n$, $\Psi(2n)\leq 4\Psi(n)$. Thus for large enough $n$, $\Phi(2n)\leq 4\Phi(n)$ and we see that if $\frac{n\ell(n)}{\Phi(n)}\underset{n\rightarrow \infty}{\longrightarrow} \infty$, then $\frac{n\ell(n)}{\Phi(2n)}\underset{n\rightarrow \infty}\longrightarrow \infty$ and $\theta_\star=\theta^{\star}=\infty$. The first statement (1) is thus established.  Since $\Phi(2n)\leq 4\Phi(n)$ and $\frac{\Phi(2n)}{2n}\geq \frac{\Phi(n)}{n}$ for large $n$, we get $\underline{\theta}^\star\geq \frac{1}{4}\overline{\theta}^{\star}$ and $\underline{\theta}^\star\leq \frac{1}{2}\overline{\theta}^{\star}$. The inequality \eqref{crudebound} readily yields $\theta^{\star}\geq \underline{\theta}^{\star}$ and we see that if  $\underset{n\rightarrow \infty}{\limsup } \ n \sum_{k=n}^{\infty}\frac{\bar{\mu}(k)}{\Phi(k)}=0$, then $\theta^{\star}\leq \overline{\theta}^{\star}$. We now establish (4). Assume that $\frac{n \ell(n)}{\Phi(n)}$ converges to $0$, since $\Phi(2n)\geq \Phi(n)$, we have that $\frac{n\ell(n)}{\Phi(2n)}\underset{n\rightarrow \infty}\longrightarrow 0$. Thus $n\sum_{k=n+1}^{2n}\frac{\bar{\mu}(k)}{\Phi(k)}\leq \frac{n\ell(n)}{\Phi(n)}$ and the left hand side of this inequality goes to $0$. We thus obtain
\[\underset{n\rightarrow \infty}{\limsup} \ n \!\!\!\sum_{k=2n+1}^{\infty}\frac{\bar{\mu}(k)}{\Phi(k)}=\underset{n\rightarrow \infty}{\limsup} \ n \sum_{k=n}^{\infty}\frac{\bar{\mu}(k)}{\Phi(k)}\]
and get an equality in \eqref{crudebound}. 
\end{proof}
As a first application of Lemma \ref{boundsigma} we get the following.
\begin{lemma}\label{suffcond2} Assume $\mu(\infty)=0$, $c_{\mathrm{k}}=0$ and $\Phi(n)\underset{n\rightarrow \infty}{\sim} dn^{1+\beta}$ for some $\beta \in (0,1)$.  If $n^{1-\beta}\bar{\mu}(n) \underset{n \rightarrow \infty}{\longrightarrow} 0$, then $\theta^{\star}=0$ and the process comes down from infinity.
\end{lemma}
\begin{proof} 
Assume $\Phi(n)\underset{n\rightarrow \infty}{\sim} dn^{\beta+1}$ for a certain constant $d>0$. Let $\epsilon>0$, by assumption there exists $N$ such that for all $k\geq N$, $\bar{\mu}(k)\leq \frac{d}{k^{1-\beta}}\frac{\epsilon}{2}.$ Hence, on the one hand, for large enough $n$, $n^{-\beta}\sum_{k=1}^{N}\bar{\mu}(k)\leq d\frac{\epsilon}{2}$.
On the other hand, since $\sum_{k=N+1}^{n}\frac{1}{k^{1-\beta}}\leq n^{\beta}$, we have that
$$n^{-\beta}\sum_{k=N+1}^{n}\bar{\mu}(k)\leq n^{-\beta} \frac{d \epsilon}{2}  \sum_{k=N+1}^{n}\frac{1}{k^{1-\beta}}\leq\frac{d \epsilon}{2},$$
thus, for large enough $n$, we have 
$\frac{n\ell(n)}{\Phi(n)}\leq d\epsilon.$
It remains to study $\underset{n\rightarrow \infty}{\limsup}\ n\sum_{k=n}^{\infty}\frac{\bar{\mu}(k)}{k^{\beta+1}}$.
One has $n\bar{\mu}(n)\sum_{k=n}^{\infty}\frac{1}{k^{\beta+1}}\underset{n\rightarrow \infty}{\sim} cn^{1-\beta}\bar{\mu}(n)$ for a certain constant $c$. By assumption,  $n^{1-\beta}\bar{\mu}(n)\rightarrow 0$ when $n\rightarrow \infty$ and since $$n\sum_{k=n}^{\infty}\frac{\bar{\mu}(k)}{k^{\beta+1}}\leq n\bar{\mu}(n)\sum_{k=n}^{\infty}\frac{1}{k^{\beta+1}}$$
by applying Lemma \ref{boundsigma}-(4), we have that $\theta=\theta^{\star}=\theta_{\star}=0$. 
\end{proof}
The following examples are easily investigated by applying Lemmas \ref{boundsigma} and \ref{suffcond2}.
\begin{example} 
\begin{itemize}
\item If $\Phi(n)\underset{n\rightarrow \infty}{\sim} dn^{\beta+1}$ with $\beta \in (0,1)$ and $\bar{\mu}(n)\underset{n\rightarrow \infty}{\sim} \frac{\lambda \log(n)^{\alpha}}{n}$ with $\alpha \in \mathbb{R}$, then $n^{1-\beta}\bar{\mu}(n)\underset{n\rightarrow \infty}{\longrightarrow} 0$ and Lemma \ref{suffcond2} ensures that the process comes down from infinity. \item If $\Phi(n)\underset{n\rightarrow \infty}{\sim} dn(\log n)^{\beta}$ with $\beta>1$ and $\bar{\mu}(n)\underset{n\rightarrow \infty}{\sim} \frac{\lambda}{n^{\alpha}}$ with $\alpha \in (0,1)$, then one can check that, for some constant $c>0$, $\frac{n\ell(n)}{\Phi(2n)}\underset{n\rightarrow \infty}{\sim} c\frac{n^{1-\alpha}}{(\log 2n)^{\beta} } \underset{n\rightarrow \infty}{\longrightarrow} \infty$ and Lemma \ref{boundsigma}-(1) ensures that the process stays infinite.
\end{itemize}
\end{example}

\subsection{Proof of Corollary \ref{FEFC2}} 
Recall the statements of Corollary \ref{FEFC2}. Set $\lambda:=\mu(\infty)\geq 0$ and $c_{\mathrm{k}}:=\Lambda(\{0\})$. We first establish (1), namely we show that if $c_{\mathrm{k}}>0$, then $\theta^\star=\theta_{\star}=\theta=\frac{2\lambda}{c_{\mathrm{k}}}\geq 0$. Recall $\bar{\mu}(k):=\mu(\{k,k+1,\ldots, \infty\})$. Let $\mu^{0}$ be the restriction of the measure to $\mathbb{N}$, for all $k\in \mathbb{N}$, $\mu^{0}(k)=\mu(k)$ and $\mu^{0}(\infty)=0$. By definition of the parameters $\theta^{\star}$ and $\theta_{\star}$ in \eqref{theta}, we get
\begin{equation}\label{thetainfefc}\theta^{\star}:=\underset{n\rightarrow \infty}{\limsup} \ \sum_{k=1}^{\infty}\frac{n\lambda}{\Phi(n+k)} +\underset{n\rightarrow \infty}{\limsup} \  \sum_{k=1}^{\infty}\frac{n\bar{\mu}^{0}(k)}{\Phi(n+k)},
\end{equation}
and similarly for $\theta_\star$ replacing $\limsup$ by $\liminf$. We study separately the two summands in \eqref{thetainfefc}. First observe that 
$\sum_{k=1}^{\infty}\frac{n\lambda }{\Phi(k+n)}=\sum_{k=n+1}^{\infty}\frac{n\lambda }{\Phi(k)}$. Recall that $\Phi(n)\underset{n\rightarrow \infty}{\sim} \Psi(n)$ with $\Psi$ given in \eqref{psi}. By \cite[Proposition 2-(i), p16]{Ber96}, one has $\frac{\Psi(n)}{n^2}\underset{n\rightarrow \infty}\longrightarrow \frac{c_{\mathrm{k}}}{2}$. If $c_{\mathrm{k}}>0$, then $\sum_{k=n}^{\infty}\frac{1}{\Phi(k)}\underset{n\rightarrow \infty }{\sim} \frac{2}{c_{\mathrm{k}}n}$ and we see plainly that 
$\underset{n\rightarrow \infty}{\limsup} \ n\lambda \sum_{k=n+1}^{\infty}\frac{1}{\Phi(k)}=2\lambda/c_{\mathrm{k}}$.
It remains to study the second summand in \eqref{thetainfefc}. We apply Lemma \ref{boundsigma}-(4) to the process with splitting measure $\mu^{0}$. Set $\ell^0(n):=\sum_{k=1}^{n}\bar{\mu}^0(k)$ for all $n\geq 1$. Since $\bar{\mu}^0(k)\underset{k\rightarrow \infty}\longrightarrow 0$ then by C\'esar\`o's theorem $\frac{\ell^0(n)}{n} \underset{n\rightarrow \infty}\longrightarrow 0$ and
$\frac{n\ell^0(n)}{\Phi(n)}\underset{n\rightarrow \infty}{\sim} \frac{2}{c_{\mathrm{k}}}\frac{\ell^{0}(n)}{n}\underset{n\rightarrow \infty}\longrightarrow 0$. Moreover, 
$n\sum_{k=n}^{\infty}\frac{\bar{\mu}^0(k)}{\Phi(k)}\leq n\bar{\mu}^0(n)\sum_{k=n}^{\infty}\frac{1}{\Phi(k)}\underset{n\rightarrow \infty}{\sim} \frac{2}{c_{\mathrm{k}}}\bar{\mu}^0(n) \underset{n\rightarrow \infty}\longrightarrow 0,$ and  $\underset{n\rightarrow \infty}{\limsup} \  \sum_{k=n+1}^{\infty}\frac{n\bar{\mu}^{0}(k)}{\Phi(k)}=0$. Finally, by Lemma \ref{boundsigma}-(4) and Lemma \ref{boundsigma}-(2), we get $\theta^{\star}=\frac{2\lambda}{c_{\mathrm{k}}}$. Similar arguments provide $\theta_{\star}=\frac{2\lambda}{c_{\mathrm{k}}}$. Recall the statement (2) of Corollary \ref{FEFC2}. Note that if $c_{\mathrm{k}}=0$ and $\lambda>0$ then $\Phi(n)/n^2\underset{n\rightarrow \infty}{\longrightarrow} 0$ and $n\ell(n)\underset{n\rightarrow \infty}{\sim} n\lambda$. Thus $\frac{n\ell(n)}{\Phi(n)}\underset{n\rightarrow \infty}{\sim} \lambda \frac{n^2}{\Phi(n)}\underset{n\rightarrow \infty}{\longrightarrow} \infty$, and Lemma \ref{boundsigma}-(1) entails that $\theta_{\star}=\infty$.

\subsection{Proof of Corollary \ref{suffcond1}} Recall the statement of Corollary \ref{suffcond1} and the definition of $\theta^{\star}$ in \eqref{theta}. Recall that we work under the assumption \eqref{CDI}. Assume  $\sum_{k=2}^{\infty}\frac{k\bar{\mu}(k)}{\Phi(k)}<\infty$. Recall \eqref{theta} and that the sequence $(k/\Phi(k),k\geq 1)$ is non-increasing. For any $n\geq 1$,
\begin{equation*}
\frac{n\bar{\mu}(k)}{\Phi(n+k)}\leq \frac{n+k}{\Phi(n+k)}\frac{n}{n+k}\bar{\mu}(k)\leq \frac{k}{\Phi(k)}\bar{\mu}(k).
\end{equation*}
Since \eqref{CDI} holds, $\Phi(n)/n\underset{n\rightarrow \infty}{\longrightarrow} \infty$. Therefore, for any $k\geq 2$, $\underset{n\rightarrow \infty}\lim \frac{n\bar{\mu}(k)}{\Phi(k+n)}=0$. Finally by Lebesgue's theorem, we have that $\theta^{\star}:= \underset{n\rightarrow \infty}{\limsup} \sum_{k=1}^{\infty}\frac{n\bar{\mu}(k)}{\Phi(k+n)}=0$.

\subsection{Proof of Proposition \ref{regularcdi}} Recall the assumptions of Proposition \ref{regularcdi}. Case (2) is a consequence of  Lemma \ref{suffcond2}. Note that in cases (1) and (3), we necessarily have $\alpha\in (0,1)$. Moreover, one has \[\bar{\mu}(n)\underset{n\rightarrow \infty}{\sim} \frac{\lambda}{n^{\alpha}},\quad\ell(n)\underset{n\rightarrow \infty}{\sim} \frac{\lambda}{(1-\alpha)}n^{1-\alpha},\ \Phi(n)\underset{n\rightarrow \infty}{\sim} dn^{\beta +1}\]
for some constant $d>0$. Clearly,
if $\alpha+\beta <1$, then \[\frac{n\ell(n)}{\Phi(n)}\underset{n\rightarrow \infty}{\sim}\frac{\lambda}{d(1-\alpha)}n^{1-(\alpha+\beta)}\underset{n\rightarrow \infty}{\longrightarrow} \infty\]
and by Lemma \ref{crudebound}-(1), one gets $\theta_\star=\infty$. We now treat the critical case, $\alpha+\beta=1$. Assume  $\bar{\mu}(k)\underset{k\rightarrow \infty}{\sim}\frac{\lambda}{k^{\alpha}}$ and $\Phi(n)\underset{n\rightarrow \infty}{\sim} dn^{2-\alpha}$. 
One returns to the definition \eqref{theta} of $\theta$. By assumption for any constants $c_1<1<c_2$, there exists $k_0$ large enough such that for all $k\geq k_0$, $$c_2\frac{\lambda}{k^{\alpha}}\geq \bar{\mu}(k)\geq c_1\frac{\lambda}{k^{\alpha}}$$ 
 and  $$c_1 d (n+k)^{2-\alpha} \leq \Phi(n+k)\leq c_2 d (n+k)^{2-\alpha} \text{ for all } n\geq 1.$$
Thus for $k\geq k_0$ and $n\geq 1$, we have
\begin{equation}\label{encadtheta} \frac{c_2}{c_1}\frac{\lambda}{d}\sum_{k=k_0}^{\infty}\frac{n}{k^{\alpha}(n+k)^{2-\alpha}}\geq \sum_{k=k_0}^{\infty}\frac{n\bar{\mu}(k)}{\Phi(n+k)}\geq \frac{c_1}{c_2}\frac{\lambda}{d}\sum_{k=k_0}^{\infty}\frac{n}{k^{\alpha}(n+k)^{2-\alpha}}.
\end{equation}
For any $n\geq 1$, $\sum_{k=1}^{k_0}\frac{n \bar{\mu}(k)}{\Phi(n+k)}\leq \frac{n}{\Phi(n)}\sum_{k=1}^{k_0}\bar{\mu}(k)$. Since $\frac{n}{\Phi(n)}\underset{n\rightarrow \infty}{\sim} \frac{d}{n^{1-\alpha}}\underset{n\rightarrow \infty}{\rightarrow} 0$, the latter bound goes to $0$ and we only need to focus on the limit as $n$ goes to $\infty$ of the series
$$\sum_{k=k_0}^{\infty}\frac{n}{k^{\alpha}(n+k)^{2-\alpha}}.$$ A comparison with an integral provides
\begin{align*}
\int_{k_0}^{\infty}\frac{n}{x^{\alpha}(n+x)^{2-\alpha}}\ddr x&\leq \sum_{k=k_0}^{\infty}\frac{n}{k^{\alpha}(n+k)^{2-\alpha}}\leq \int_{k_0-1}^{\infty}\frac{n}{x^{\alpha}(n+x)^{2-\alpha}}\ddr x.
\end{align*}
By factorizing $n$ and doing the change of variable $u=\frac{x}{n}$, we get
\begin{align*}
\int_{k_0/n}^{\infty}\frac{1}{u^{\alpha}(1+u)^{2-\alpha}}\ddr u&\leq \sum_{k=2}^{\infty}\frac{n}{k^{\alpha}(n+k)^{2-\alpha}}\leq \int_{(k_0-1)/n}^{\infty}\frac{1}{u^{\alpha}(1+u)^{2-\alpha}}\ddr u.
\end{align*}
Letting $n$ to $\infty$, provides
\begin{align*}
\underset{n\rightarrow \infty}\lim \sum_{k=k_0}^{\infty}\frac{n}{k^{\alpha}(n+k)^{2-\alpha}}&=\int_{0}^{\infty}\frac{\ddr u}{u^{\alpha}(1+u)^{2-\alpha}}
=\left[\frac{1}{1-\alpha}\left(\frac{u}{u+1}\right)^{1-\alpha}\right]_{0}^{\infty}=\frac{1}{1-\alpha}.
\end{align*} 
We deduce from \eqref{encadtheta} that \[\theta^{\star}, \theta_{\star}\in \left[\frac{c_1}{c_2}\frac{\lambda}{d(1-\alpha)},\frac{c_2}{c_1}\frac{\lambda}{d(1-\alpha)}\right].\] Since $c_1$ and $c_2$ can be chosen arbitrarily close to $1$, we get $\theta^{\star}=\theta_{\star}=\frac{\lambda}{d(1-\alpha)}$.

\subsection{Proof of Proposition \ref{slowcdi}} 
We now consider some simple EFC processes with ``slow" coalescence.  Recall the assumptions. As $n$ goes to $\infty$,  $\Phi(n)\sim dn(\log n)^{\beta}$ with $\beta>1$ and $\bar{\mu}(n)\sim \lambda (\log n)^{\alpha}/n$ with $\alpha\in \mathbb{R}$. 

If $\alpha<0$ then $\beta-\alpha>1$ and since $\frac{n\bar{\mu}(n)}{\Phi(n)}\underset{n\rightarrow \infty}{\sim} \frac{\lambda}{d}\frac{1}{n (\log n)^{\beta-\alpha}}$, we have that $\sum_{k=2}^{\infty}\frac{k\bar{\mu}(k)}{\Phi(k)}<\infty$. Applying Corollary \ref{suffcond1} entails that $\theta^{\star}=0$. 

We now focus on the case $\alpha>0$ and to simplify the calculations, we treat it with the assumption $\bar{\mu}(n)=\frac{\lambda (\log n)^{\alpha}}{n}$ for any $n\geq 2$. Plainly since $\beta>1$, then \eqref{CDI} holds and one can apply Theorem \ref{mainthm}. We now compute $\theta$. 
A  comparison with integrals provides 
\[n\int_{1}^{n}\frac{(\log x)^{\alpha}}{x+1}\ddr x\leq n \ell(n)\leq n\int_{1}^{n}\frac{(\log(x+1))^{\alpha}}{x}\ddr x.\]
Both integrands are equivalent to $(\log x)^{\alpha}/x$ as $x$ goes to $\infty$ and we get
$n\ell(n)\underset{n\rightarrow \infty}{\sim} \frac{\lambda}{\alpha +1} n(\log n)^{\alpha+1}$. One checks
\begin{equation}\label{equiv} \frac{n\ell(n)}{\Phi(n)}\underset{n\rightarrow \infty}{\sim} \frac{\lambda}{d(\alpha +1)}(\log n)^{1-(\beta-\alpha)}.\end{equation}
We now apply Lemma \ref{boundsigma}.
\begin{enumerate}
\item If $\beta-\alpha<1$, then
$\frac{n\ell(n)}{\Phi(n)}\underset{n\rightarrow \infty}{\longrightarrow} \infty$ and by Lemma \ref{boundsigma}-(1), $\theta_{\star}=\infty$.
\item If $\beta-\alpha\geq 1$, then  notice first that for some constant $C>0$ and large enough $n$,  
\[n\sum_{k=n}^{\infty}\frac{\bar{\mu}(k)}{\Phi(k)}\leq n\sum_{k=n}^{\infty}\frac{C}{k^{2}(\log k)^{\beta-\alpha}}\leq \frac{C}{\log n}.\]
Therefore, $n\sum_{k=n}^{\infty}\frac{\bar{\mu}(k)}{\Phi(k)}\underset{n\rightarrow \infty}{\longrightarrow} 0$. If $\beta-\alpha>1$, then we see from \eqref{equiv}, that $\overline{\theta}^{\star}:=\underset{n\rightarrow \infty}{\limsup} \frac{n\ell(n)}{\Phi(n)}=0$ and by Lemma \ref{boundsigma}-(2), $\theta^{\star}=0$. In the case of equality, $\beta-\alpha=1$, Lemma \ref{boundsigma}-(2) provides $\theta^\star\leq \frac{\lambda}{d(1+\alpha)}$. By definition of $\theta_\star$, for any $r>0$,
\begin{equation} \label{lowerbound} \theta_\star\geq \underset{n\rightarrow \infty}{\liminf} \sum_{k=1}^{\left\lfloor \frac{n}{r}\right \rfloor}\frac{n\bar{\mu}(k)}{\Phi(n+k)}\geq \underset{n\rightarrow \infty}{\liminf} \frac{n}{\Phi(n(1+1/r))}\sum_{k=1}^{\left \lfloor \frac{n}{r}\right \rfloor }\bar{\mu}(k).\end{equation}
Since $\sum_{k=1}^{\lfloor \frac{n}{r}\rfloor}\bar{\mu}(k)\underset{n\rightarrow \infty}{\sim}\frac{\lambda}{1+\alpha}\log\left(\lfloor \frac{r}{n}\rfloor \right)^{\alpha+1}$ and $\Phi(n)\underset{n\rightarrow \infty}{\sim} dn(\log n)^{\alpha+1}$, the right-hand side in \eqref{lowerbound} equals $\frac{\lambda}{d(1+\alpha)}\frac{1}{1+1/r}$ and therefore $\theta_\star\geq  \frac{\lambda}{d(1+\alpha)}\frac{1}{1+1/r}$.  Since $r$ is arbitrarily large, $\theta_\star\geq \frac{\lambda}{d(1+\alpha)}$ and \[\theta_\star=\theta^{\star}=\frac{\lambda}{d(1+\alpha)}.\]
\end{enumerate}
The case where only the equivalence  $\bar{\mu}(n)\underset{n\rightarrow \infty}{\sim}\frac{\lambda (\log n)^{\alpha}}{n}$ holds, follows from an adaptation of the previous calculations. This ends the proof of Proposition \ref{slowcdi}.\\ 

%
We conclude this article with a few comments. When there are no sudden fragmentations into infinitely many blocks, i.e. $\mu(\infty)=0$, the question whether $(\#\Pi(t),t\geq 0)$ starting from a finite state $n$ can reach $\infty$ has not been addressed.  In particular, we mention that the condition $\theta>0$  does not imply the explosion in general. Explosion requires a study in its own right by designing other taylor-made criteria. This is investigated in the work of Foucart and Zhou \cite{explosion}. We mention also the work \cite{FoucartZhouWF}  where certain Markov processes in duality with simple EFCs, called Wright-Fisher processes with selection, are studied. Other properties of the block-counting process, such as its Feller property, are stated in this latter work.


\subsection*{Acknowledgements} 
I am grateful to Bastien Mallein for many insightful discussions. I would also like to thank Martin M\"ohle and Xiaowen Zhou to whom I spoke about this problem in 2014 and 2019 respectively. 
This research has been supported  by LABEX MME-DII (ANR11-LBX-0023-01). 

\end{document}